\documentclass[12pt,eqno]{article}
\usepackage{amsmath}
\usepackage{amsfonts}
\usepackage{mathrsfs}
\usepackage{amssymb,amsmath,latexsym}

 \usepackage[dvips]{color}
\usepackage{amsmath,amsthm}     
\usepackage{amssymb}            
\usepackage{euscript}           
\usepackage{graphicx}
\usepackage{tikz}
\usepackage{curves,enumerate}

\usepackage[colorlinks]{hyperref}

\newtheorem{theorem}{Theorem}[section]

\newtheorem{lemma}[theorem]{Lemma}
\newtheorem{remark}{Remark}
\newtheorem{example}{Example}[section]
\newtheorem{corollary}[theorem]{Corollary}
\newtheorem{definition}{Definition}[section]
\newtheorem{proposition}{Proposition}[section]
\newtheorem{hyp}{Hypothesis}

\def\lab#1{\mbox{$\label{#1}$}}

\def\ga{\gamma}                         

\def\ga{\gamma}
\def\ra{\rightarrow}                    
\def\proof{\noindent {\bf Proof: }\ }

\def\lab{\label}
\def\beq{\begin{equation}}               
\def\eeq{\end{equation}}                 
\def\bea{\begin{eqnarray}}             
\def\eea{\end{eqnarray}}               
\def\be*{\begin{eqnarray*}}             
\def\ee*{\end{eqnarray*}}               
\def\ba{\begin{array}}                  
\def\ea{\end{array}}                    

\def\;{\vspace{3mm} \\ }

\def\R{\mathbb{R}}

\def\E{\mathbb{E} }
\def\P{\mathbb{P} }

\def\~{\widetilde}

\def\beqlb{\begin{eqnarray}} \def\eeqlb{\end{eqnarray}}
\def\beqnn{\begin{eqnarray*}} \def\eeqnn{\end{eqnarray*}}

\def\<{\langle}  \def\>{\rangle}

\def\bB{\mbox{\boldmath $B$}} 
\def\A{{\cal A}} \def\B{{\cal B}}

\def\ra{\rightarrow}                    

\def\bde{\begin{definition}}
\def\ede{\end{definition}}

\def\bth{\begin{theorem}}
\def\eth{\end{theorem}}
\def\bprop{\begin{proposition}}
\def\eprop{\end{proposition}}
\def\ble{\begin{lemma}}
\def\ele{\end{lemma}}
\def\bcor{\begin{corollary}}
\def\ecor{\end{corollary}}

\begin{document}

\title{\Large {\bf Joint H{\"o}lder continuity of local time for a class of interacting branching measure-valued diffusions}
\thanks{Partial funding in support of this work was provided by the Natural
Sciences and Engineering Research Council of Canada (NSERC) and the Department of Mathematics at the University of Oregon.}}
\author{\bf D. A. Dawson\\Carleton University \and \bf J. Vaillancourt\\HEC Montr\'eal 
\and \bf H. Wang\\University of Oregon}

\date{}
\maketitle

{\narrower{\narrower

\centerline{\bf Abstract}

\bigskip

Using a Tanaka representation of the local time for a class of superprocesses with dependent spatial motion, 
as well as sharp estimates from the theory of uniformly parabolic partial differential equations, the joint H\"{o}lder continuity in time and space of said local times is obtained in two and three dimensional Euclidean space. 
\bigskip

\noindent{\bf 2010 Mathematics Subject Classification}: Primary 60J68, 60J80; 
Secondary 60H15, 60K35, 60K37


\noindent{\bf Keywords and Phrases}: Measure-valued diffusions, stochastic partial differential equations, superprocesses, branching processes, joint H{\"o}lder continuity, local time, Tanaka formula.

\par}\par}


 \baselineskip=6.0mm


\section{Introduction}\label{sec:intro}

\setcounter{equation}{0}

Given a stochastic process $\{\mu_t:t\ge0\}$ built on a probability space $(\Omega, {\cal F}, \P)$ 
and valued in the space $M(\R^d)$ of all positive Radon measures on $\R^d$, and writing 
${\cal B}(\R^d \times[0,\infty))$ for the Borel $\sigma$-algebra over $\R^d \times[0,\infty)$,
a ${\cal B}(\R^d \times[0,\infty)) \times {\cal F}$-measurable function 
  $\Lambda^{x}_t(w):(\R^d \times [0, \infty) \times \Omega) \rightarrow [0, \infty)$
 is called a local time of $\{\mu_t\}$ if,  
 for any continuous function with compact support $\phi \in C_c(\R^d)$, there holds 
 $\P$-almost surely at each time $t\ge0$ (with finite integrals on both sides) 
 \beqlb
 \label{LT}
 \int_{\R^d}\phi(x) \Lambda_t^x dx = \int_0^t \<\phi, \mu_s\>ds. 
 \eeqlb
Here and henceforth we write $\<\phi,\mu\>= \int_{\R^d} \phi(x) \mu(dx)$, for any $\mu \in M(\R^d)$ 
and any $\mu$-integrable $\phi$. 
If in addition $\Lambda_t^x$ is integrable with respect to Lebesgue measure $\lambda_0$ on $\R^d$, 
for each fixed $t$, then it is the Radon-Nikodym derivative (in space) with respect to $\lambda_0$ 
of the occupation time process $\int_0^t \mu_s ds$, itself a time-averaged measure-valued process 
with (hopefully) more regular paths than $\mu_t$ itself. 

For example, Super-Brownian motion $\{\mu_t\}$ has a density
 $f_t=d\mu_t/d\lambda_0$ when $d=1$ (Konno and Shiga \cite{KonnoShiga88})
 but not in the cases $d\ge2$ (see Dawson and Hochberg \cite{DawsonHochberg79} and 
 Perkins \cite{Perkins88} \cite{Perkins89}). 
 When $d=1$ the choice $\Lambda_t^x= \int_0^t f_s(x)ds$ is immediate but when $d\ge2$ 
 the fractal dimension of support for $\mu_t$ renders the proof of existence of $\Lambda_t^x$ 
 a lot more involved. The sharpest estimates for the closed support of Super-Brownian motion are found in 
 Dawson and Perkins \cite{DawsonPerkins91} when $d\ge3$ and
 Le Gall and Perkins \cite{LeGallPerkins95} when $d=2$. 

The existence and the joint space-time continuity of paths for the local time of Super-Brownian motion, 
when $d\le3$, were first obtained by Iscoe \cite{Iscoe86b} and Sugitani \cite{Sugitani89}. 
Their results were variously sharpened and generalized, first by Adler and Lewin \cite{AdlerLewin92} for super stable 
processes and Krone \cite{Krone93} for superdiffusions; then by many others, most notably 
Ethier and Krone \cite{EthierKrone95} for some related Fleming-Viot processes with diffusive mutations;  
L\'opez-Mimbela and Villa \cite{Lopez-MimbelaVilla04}, who streamlined and unified the various definitions of the local time and clarified their interrelations in the above cases; Li and Xiong \cite{LiXiong07}, who offered an alternative (trajectorial) definition of the local time when the superprocess is degenerate, that is, a purely atomic measure-valued process, and proved its joint H{\"o}lder continuity, as well as scaling limit theorems, using a representation in terms of stochastic integrals with respect to the excursions of an underlying Poisson random measure. 

More recently, Dawson et al. \cite{DVW19} proved the existence of the local time for the (larger still) class of superprocesses with dependent spatial motion (SDSM) first introduced in 
Wang \cite{Wang97}, \cite{Wang98}. 
This class, which contains Super-Brownian motion and some superdiffusions as special cases, describes the movement of a cloud of infinitely many particles, each subject to critical branching and diffusing according to second order partial differential operators with spatially dependent coefficients, among a random environment prescribed by a Brownian sheet. Following Hochberg \cite{Hochberg91} in his description of Super-Brownian motion, the resulting $M(\R^d)$-valued process can been colourfully described as the high density limit for an instantaneous smearing of the diffusive trajectories performed by approximating clusters of finitely many newborn particles, further subjected to the motion of an ambiant random medium. 

Our main result (Theorem \ref{lt_th1}) is the joint H{\"o}lder continuity for the local time $\Lambda_t^x$ of SDSM, 
in time and space $(t,x)$, when $d\le3$. It is formulated next, in Section \ref{sec:main}, 
together with all the notation and preliminary statements required for this purpose.  
The proof of our main result is in Section \ref{sec:Tanaka} and uses: a duality argument 
for the evaluation of the moments of $\{\mu_t\}$, inspired by Sugitani \cite{Sugitani89}, 
Krone \cite{Krone93} and Ethier and Krone \cite{EthierKrone95}; 
sharp inequalities due to Aronson \cite{Aronson68} and Lady\u{z}enskaja et al. \cite{LSU68} 
for the fundamental solutions of those uniformly parabolic partial differential equations associated with the 
Markov semigroup for SDSM; estimates for the Green function of this semigroup; and a Tanaka formula from 
Dawson et al. \cite{DVW19}, which involves stochastic integrals with singular integrands.  
Some incidental proofs are relegated to Section \ref{app:ProofsOfLemmas}. 


\section{Conditions, estimates and main results} \label{sec:main}

\subsection{Basic notation}

For any topologically complete and separable metric space $S$ (hereafter, a Polish space),
${\cal B}(S)$ denotes its Borel $\sigma$-field, $B(S)$ the Banach space
of real-valued bounded Borel measurable functions on $S$ with the supremum norm $\|\cdot\|_{\infty}$
and $C(S)$ the space of real-valued continuous functions on $S$. 
Subscripts $b$ or $c$ on a space of functions refers to its subspace of bounded
or compactly supported functions, respectively, as in $C_b(S)$ and $C_c(S)$. 
As usual $C([0, \infty), S)$ denotes the Polish space of continuous trajectories into $S$ 
with the topology of uniform convergence on compact time sets. 
$S^m$ denotes the $m$-fold product of $S$. 

When $S=\R^d$, subscript $0$ indicates those functions vanishing at $\infty$ while superscript $k\ge1$ means 
continuous derivatives up to and including order $k$ (possibly infinite). We shall make use mostly of 
$C_0(\R^d)$, those bounded continuous functions vanishing at $\infty$, its subspaces 
$C_c^2(\R^d)\subset C_0^2(\R^d)$ and the chain 
$C_c^\infty(\R^d)\subset C_b^\infty(\R^d)\subset C_b^2(\R^d)\subset C^2(\R^d)$. 
The main set of functions of interest here is 
\begin{equation*}
K_a(\R^d)=\{\phi: \phi= h + \beta I_a, \beta \in \R, h \in C_c^{\infty}(\R^d)\},
\end{equation*} 
defined for any real number $a\ge0$ with $I_a(x)= (1 + |x|^2)^{(-a/2)}$ 
 and $|x|^2 = \sum_{i=1}^dx_i^2$.

We denote by ${\rm Lip}(\R^d)$ the space 
of Lipschitz functions on $\R^d$, that is, $\phi \in {\rm Lip}(\R^d)$ 
if there is a constant $M>0$ such that $|\phi(x)-\phi(y)|\leq M |x-y|$ 
for every $x, y \in \R^d$. Its subset of bounded functions is written ${\rm Lip}_b(\R^d)$. 

We will also need $C_b^{1,2}([0, t] \times (\R^d)^m)$, the space of 
bounded continuous functions with all derivatives bounded, up to and including order $1$ in the time variable up to time $t$
and order $2$ in the $md$ space variables, including mixed derivatives of that order. 
When no ambiguity is present we also write the partial derivatives 
(of functions and distributions) in abridged form
\begin{equation*}
\partial_p=\frac{\partial}{\partial x_p}, \quad\quad 
\partial_p^2=\frac{\partial^2}{\partial x_p^2} \quad \mbox{and}\quad 
\partial_p\partial_q=\frac{\partial}{\partial x_p}\frac{\partial}{\partial x_q} \quad \mbox{and so on.}
\end{equation*}

Given any positive Radon measure $\mu\in M(\R^d)$ and any $p\in[1,\infty)$, 
 we write $L^p(\mu)$ for the Banach space of real-valued Borel measurable functions on $\R^d$, 
 with finite norm $\|\phi\|_{\mu,p}:= \{\int_{\R^d}|\phi(x)|^pd\mu(x)\}^{1/p} < \infty$. 
  
Let $M(\R^d)$ be the space of all positive Radon measures on $\R^d$ and 
$M_0(\R^d)$, its subspace of finite positive Radon measures. 
For any real number $a\ge0$, define the main set of measures of interest here as  
\begin{equation*}
M_a(\R^d)= \{\mu \in M(\R^d): \<I_a,\mu\>= \int_{\R^d}I_a(x) \mu(dx) < \infty \}.
\end{equation*}  
The topology $\tau_a$ of $M_a(\R^d)$ is defined in the following way: 
$\mu_n \in M_a(\R^d) \Rightarrow \mu \in M_a(\R^d)$ as $n \ra \infty$, iff $\lim_{n \ra \infty}\<\phi , \mu_n \> =
\<\phi , \mu \>$ holds for every $\phi \in K_a({\R^d})$. Then, $(M_a(\R^d), \tau_a)$ is a Polish space 
(see Iscoe \cite{Iscoe86a} and Konno and Shiga \cite{KonnoShiga88}). 
For instance, the Lebesgue measure $\lambda_0$ on $\R^d$ belongs to $M_a(\R^d)$ for any $a >d$. 
Furthermore, both $dx=\lambda_0(dx)$ are used indifferently when calculating Lebesgue integrals. 

The short form $\mu^{m}= \mu\otimes\ldots\otimes\mu$
denotes the $m$-fold product measure of $\mu\in M_a(\R^d)$ by itself
and we write ${\cal I}_{a, m}$ for the product  
${\cal I}_{a, m}(x)=I_a(x_1)\cdot\ldots\cdot I_a(x_m)$,  
keeping in mind that 
${\cal I}_{a, m}^{-1} f(x)={\cal I}_{a, m}^{-1}(x)\cdot f(x)$ 
means the product, not the composition of functions. 


\subsection{Sufficient conditions}

The following basic assumptions are valid throughout this paper.

\begin{hyp}\label{hyp:basicassumpFilter}
Let $(\Omega, {\cal F}, \{{\cal F}_t\}_{t \geq 0}, \P)$ be a filtered probability space with a right continuous filtration
$\{{\cal F}_t\}_{t \geq 0}$, satisfying the usual hypotheses and upon which all our processes are built, notably a  
Brownian sheet $W$ on $\R^d$ and a countable family $\{B_{k}, k\geq 1\}$ of independent, 
$\R^d$-valued, standard Brownian motions written $B_{k}=(B_{k1}, \cdots, B_{kd})$. 
The family $\{B_{k}, k\geq 1\}$ is assumed independent of $W$. 
\end{hyp}

Recall that a Brownian sheet on $\R^d$ is an $\R^1$-valued random set function $W$ on the 
Borel $\sigma$-field ${\cal B}(\R^d \times \R_+)$ defined on 
$(\Omega, {\cal F}, \{{\cal F}_t\}_{t \geq 0}, \P)$ such that both of the following statements hold: 
for every $A \in {\cal B}(\R^d)$ having finite Lebesgue measure $\lambda_0(A)$, 
the process $W(A\times [0, t])$ is a square-integrable $\{ {\cal F}_t\}$-martingale;
and for every pair $A_i\in {\cal B}(\R^d \times \R_+)$, $i=1, 2$, having finite
Lebesgue measure with $A_1\cap A_2=\emptyset$, the random variables $W(A_1)$ and
$W(A_2)$ are independent, Gaussian random variables with mean zero, 
respective variance $\lambda_0(A_i)$ and
$ W(A_1\cup A_2) = W(A_1)+W(A_2)$ holds $\P$-almost surely  
(see Walsh \cite{Walsh86} or Perkins \cite{Perkins02}). 

The mathematical description of the diffusive motion of SDSM requires the following 
second order partial differential operators: for all $f\in C_b^{2}({(\R^d)^m})$, 
 \beqlb\label{eqn:Gn}
{G}_{m} f(\bar x):=
   \frac{1}{2} \sum_{i,j=1}^{m} \hspace{1mm}
 \sum_{p, q=1}^{d} \Gamma_{pq}^{ij}(\bar x)
  \frac{\partial^2}{\partial x_{ip} \partial x_{jq} } f(\bar x)
 \eeqlb
 where $\bar x=(x_1, \cdots, x_m)\in (\R^d)^m$ has components 
 $x_i = (x_{i1}, \cdots, x_{id}) \in \R^d$ for $1\leq i \leq m$ and 
 $\Gamma_{pq}^{ij}$ is defined by 
   \bea \label{gammaij} \hspace*{8mm}\Gamma_{pq}^{ij}(\bar x) := \left\{
\begin{array}{ll}
                                 (a_{pq}(x_{i})+\rho_{pq}(0))
                                 \quad & \mbox{if
                                 $i=j$,} \\
                                 \rho_{pq}(x_{i}- x_{j}) \quad & \mbox {if $i\neq j$},
                                 \end{array}
                                 \right.
\eea
where, for $x, y \in \R^d$ and $p,q=1,\ldots,d$, the local (or individual) diffusion coefficient is
 \beqlb \lab{a_pq}
  a_{pq}(x):=\sum_{r=1}^{d}c_{pr}(x)c_{qr}(x)
 \eeqlb
 and the global (or common) interactive diffusion coefficient is 
 \beqlb \lab{rho_pq}
 \rho_{pq}(x-y):=\int_{\R^d}h_{p}(u-x)h_{q}(u-y)du.
\eeqlb 
 
 \begin{hyp}\label{hyp:basicassumpElliptic} 
 The vector $h=(h_1, \cdots,h_d)$ satisfies $h_i \in  L^1(\R^d)\cap {\rm Lip}_b(\R^d)$    
and the $d\times d$ matrix $c=(c_{ij})$ satisfies $c_{ij} \in{\rm Lip}_b(\R^d)$, for every $i,j=1, \cdots,d$. 
 For every $m\ge1$, the $dm\times dm$ diffusion matrices $(\Gamma_{pq}^{ij})_{1 \leq i,j \leq m; 1 \leq p,q \leq d}$ 
 of real-valued functions defined by $(\ref{gammaij})$ are strictly positive definite everywhere on $(\R^d)^m$, that is, 
 there exists two positive constants $\lambda^*$ and $\Lambda^*$ such that for any 
$\xi = (\xi^{(1)}, \cdots, \xi^{(m)}) \in (\R^d)^m$ we have 
\beqlb \lab{UEC} 
 0 <  \lambda^* |\xi|^2\leq \sum_{i,j=1}^m \sum_{p,q=1}^d \Gamma_{pq}^{ij}(\cdot) \xi_p^{(i)} \xi_q^{(j)} 
 \leq \Lambda^*|\xi|^2 < \infty.
\eeqlb
\end{hyp}
 
 Hypothesis \ref{hyp:basicassumpElliptic} ensures that operator ${G}_{m}$ is not only uniformly elliptic, 
 but also that it generates a Feller semigroup $\{P^m_t: t\ge 0\}$ mapping each of 
$B({(\R^d)^m})$, $C_b({(\R^d)^m})$ and $C_0({(\R^d)^m})$ into itself, for each $t>0$ 
(see Ethier and Kurtz \cite[Chapter 8]{EthierKurtz86}). It can be written as  
 \beqlb\label{eqn:Semigroup}
 P_t^m f (\cdot)=\int_{(\R^d)^m} f(y) q_t^{m}(\cdot,y)dy
 \eeqlb
when $t > 0$, for every $f\in C_0({(\R^d)^m})$, with a transition probability density 
$q_t^{m}(x,y) > 0$ which is jointly continuous 
in $(t,x,y) \in (0, \infty) \times (\R^d)^m \times (\R^d)^m$ 
everywhere and such that $P_t^m f (x)$, as a function of $(t,x)$, 
belongs to $\cup_{t\ge0}C_b^{1,2}((0, t] \times (\R^d)^m)$, 
for every choice of $f\in C_0({(\R^d)^m})$ 
(see Stroock and Varadhan \cite[Chapter 3]{StroockVaradhan79}). 

The rest of the assumptions collectively constrain the family of initial measures. 
When the initial state of SDSM is a finite 
measure, the total mass process of SDSM is equivalent to a one-dimensional continuous state critical 
branching process and extinction occurs almost surely (see Wang \cite{Wang97}, \cite{Wang98}). 
Unbounded measures as initial states give rise to 
alternative phenomena, provided that the potential explosion of mass near the initial time $t=0$ is curtailed 
in order to allow the local time of SDSM to exist at all (see Dawson et al. \cite{DVW19}). 
This explosion can occur in either situation, as exemplified next.

\begin{example}\lab{explosionat0} Let $\varphi_s$ be the transition density of a Brownian motion particle on $\R^d$
\beqlb \lab{gaussiankernel}
   \varphi_s(y):= \frac{1}{(2 \pi s)^{d/2}} \exp{\{- \frac{|y|^2}{2s}\} }.
  \eeqlb
When the initial measure is $\delta$, the Dirac measure which puts mass $1$ 
at the origin $0\in\R^d$, we get, for all $t >0$ and $d\ge1$,
$ \sup_{0 < s \leq t} \<\varphi_s(y-\cdot), \delta\> < \infty$ if and only if $y \neq 0$.
\end{example}

Under less restrictive conditions than Hypothesis \ref{hyp:basicassumpElliptic}, 
Aronson \cite[Theorem 10]{Aronson68} proved that, for the transition probability density 
$q^{1}$ from (\ref{eqn:Semigroup}), there exist four positive constants 
$a^*$, $b$, $c$ and $A^*$ such that
\bea \label{Aronsonbounds} 
   a^* \cdot  \varphi_{bs} (y-x) \leq q_s^{1}(x,y) \leq A^* \cdot \varphi_{cs}(y-x)
\eea
holds for any $x, y \in \R^d$ and $s > 0$. 
Example \ref{explosionat0} thus remains valid in the uniformly elliptic case as well:  
an initial measure with an atom makes the existence of the  
local time for SDSM questionable in general, hence the need for an additional restriction. 

\begin{hyp}\label{hyp:basicassumpGaussUniform}
For all $t > 0$, the initial measure $\mu_0\in\cup_{b\ge0}M_b(\R^d)$ verifies
\bea \label{nonBAUnif} 
\Upsilon_t(q^{1}, \mu_0):=\sup_{y\in \R^d} \sup_{0 < s \leq t} \<q_s^{1}(\cdot,y), \mu_0\> < \infty.
\eea
\end{hyp} 

 \noindent  {\bf Remark:} 
 Any measure $\mu_0$ which is absolutely continuous with respect to 
$\lambda_0$ and has a Radon-Nikodym derivative 
which is either bounded or finitely $\lambda_0$-integrable, satisfies 
Hypothesis \ref{hyp:basicassumpGaussUniform}, notably the measures $I_a(x)dx$ 
for all choices of $a\ge0$. Combined with Aronson's inequalities (\ref{Aronsonbounds}), 
the requirement of finiteness in (\ref{nonBAUnif}) is equivalent to its special case 
$\Upsilon_t(\varphi, \mu_0)< \infty$ obtained by setting 
$q_s^{1}(\cdot,y)=\varphi_s(y-\cdot)$ from (\ref{gaussiankernel}). 

The proof of joint continuity of the local time most likely requires such a uniform bound, since our efforts 
at localizing this to a neighbourhood of the origin, succeed at the cost of a strengthening of the observation that, 
for every $w\in\R^d$, the translate $\mu_0(\cdot-w)$ 
of any measure $\mu_0 \in M_a(\R^d)$ is also in $M_a(\R^d)$. This is immediately seen though the inequality 
$I_a(x+w) \leq 2^{(a/2)} I_a^{-1}(w)I_a(x)$, valid for every $x,w\in \R^d$. 

\begin{hyp}\label{hyp:basicassumpUniformInteg} 
The initial measure $\mu_0\in\cup_{b\ge0}M_b(\R^d)$ verifies
\[
 \sup_{w\in \R^d}  \<I_a(\cdot+w), \mu_0 \> < \infty 
\] 
for some $a\ge0$. 
\end{hyp} 
This is true for any finite measure $\mu_0 \in M_0(\R^d)\subset M_a(\R^d)$ and every $a\ge0$, 
hence any Radon measure with compact support and any measure with a 
finitely $\lambda_0$-integrable Radon-Nikodym derivative with respect to $\lambda_0$; 
as well as any measure with a bounded Radon-Nykodym derivative with respect to 
$\lambda_0$, in this last case provided $a >d$. 

Finally, observe that when Hypotheses \ref{hyp:basicassumpElliptic}  
and \ref{hyp:basicassumpUniformInteg} are combined, 
Hypothesis \ref{hyp:basicassumpGaussUniform} is equivalent to the 
weaker requirement of the existence of an  
$\epsilon > 0$ such that 
\bea \label{nonBA} 
\sup_{y\in \R^d} \sup_{0 < s \leq \epsilon}\<\varphi_s(y-\cdot) , \mu_0\> < \infty
\eeqlb
holds. 
The proof of this last statement is in Subsection \ref{app:pf_nonBAUnif}.  


 \subsection{Sharp estimates}

We need some important properties of fundamental solutions to uniformly parabolic partial differential equations, 
specifically here, the Kolmogorov backward equation $\mathbb{L}u = 0$, recalling that operator
$\mathbb{L} = {G}_{1}-\partial_t$ is uniformly parabolic whenever ${G}_{1}$ is uniformly elliptic. 
 
\ble \lab{le_new0}
Under Hypothesis $\ref{hyp:basicassumpElliptic}$ and for each choice of $T > 0$ and $d\ge1$, 
there are positive constants $a_1$ and $a_2$ such that,
for all nonnegative integers $r$ and $s$ verifying $0\le 2r+s\le2$,
\beqlb \lab{LSU}
   \left|\frac{\partial^{r}}{\partial t}\frac{\partial^{s}}{\partial y_{p}} q_t^{1}(x,y)\right|
   \leq \frac{a_1}{t^{(d+2r+s)/2}}\exp{\left\{- a_2\left( \frac{|y-x|^{2}}{t} \right)\right\}}
  \eeqlb
 holds everywhere in $(t,x,y)\in(0,T) \times \R^d \times \R^d$ with $y=(y_1,\ldots,y_{d})$. 
 There also exists a unique fundamental solution $\Gamma(x,t; \xi, \tau)$ to $\mathbb{L}u = 0$,  
 started at position $\xi\in\R^d$ at initial time $\tau\in[0,T)$ and evaluated at new position $x\in\R^d$ 
 after spending time $t-\tau > 0$ to reach $t\in(\tau,T]$.
Moreover, there exist two constants $c > 0$ and $c_0 > 0$, such that, 
for all nonnegative integers $r,s_1,s_2,\ldots,s_d$ verifying 
$0\le l=2r + s_1+s_2+ \ldots + s_d \le 2$ and writing 
$\partial^{l}=\partial^{r}_t\partial^{s_1}_{x_1}\partial^{s_2}_{x_2}\cdots\partial^{s_d}_{x_d}$  
with $\partial^{0}$ for the identity, there holds, for every choice of $\alpha\in(0,1)$,
 \beqlb \lab{formulaI}
 & & |\partial^{l}\Gamma(x, \xi; t, \tau) -\partial^{l}\Gamma(x, \xi^{'}; t, \tau^{'})|  \nonumber  \\
\leq & & c( |\xi-\xi^{'}|^{\alpha} + |\tau - \tau^{'}|^{\alpha/2}) 
\bigg[ (t-\tau)^{- (d+l)/2} \exp{ \{- c_0 \frac{|x - \xi|^2}{t - \tau} \} } \nonumber  \\
 & & \hspace{5cm} + (t-\tau^{'})^{- (d+l)/2} \exp{\{-c_0 \frac{|x - \xi^{'}|^2}{t - \tau^{'}}\}}  \bigg]. 
 \eeqlb
 \ele
  \proof The upper bound in (\ref{LSU}) is a consequence of equation (13.1) of 
  Lady\u{z}enskaja et al. (\cite{LSU68} p.376). The rest of this lemma is a special case of 
  Theorem 3.5 from Garroni and Menaldi \cite[Chapter 5, Section 3]{GarroniMenaldi92}.
  Alternatively one can prove this directly using a sequence of estimates from 
  Lady\u{z}enskaja et al. \cite[Chapter IV, Sections 11,12,13]{LSU68}. 
  \qed
    
Prototypes of these bounds go back to Dressel \cite{Dressel40}, \cite{Dressel46} 
and Aronson \cite{Aronson68}.
  
   
\subsection{SPDE and dual process}

Let us now describe the SDSM class itself. 

Under Hypotheses \ref{hyp:basicassumpFilter} and \ref{hyp:basicassumpElliptic}, 
 Dawson et al. \cite{DVW19} characterize SDSM by way of a well-posed martingale problem, the unique 
solution of which is the law of a measure-valued diffusion (Markov process with continuous paths) 
$\{\mu_t:t\ge0\}$ which satisfies
 \beqlb \label{SPDEb}
  \<\phi,\mu_{t}\>  -  \<\phi,\mu_{0}\> 
  = {X}_{t}(\phi)   + M^\mu_t(\phi)  + \int_{0}^{t} \left\<  G_1\phi, \, {\mu}_{s} \right\> ds 
 \eeqlb
for every $t>0$, $\phi \in K_a(\R^d)$ and $\mu_{0}\in M_a(\R^d)$, where both 
$$
{X}_{t}(\phi)  := \sum_{p=1}^{d} \int_{0}^{t}\int_{\R^d}
\left\<h_p(y-\cdot) \partial_p \phi(\cdot), {\mu}_{s} \right\>
{W}(dy,ds)
$$
and 
$$ M^\mu_t(\phi) :=\int_0^t \int_{\R^d} \phi (y) M^\mu(ds, dy)
$$
are continuous square-integrable $\{{\cal F}_t\}$-martingales, mutually orthogonal for 
every choice of $\phi \in K_a(\R^d)$ and driven respectively by a Brownian sheet $W$ 
and a square-integrable martingale measure $M^\mu$ (generated by process $\{\mu_t:t\ge0\}$) with quadratic variation 
 \beqnn
  \<{M^\mu}(\phi) \>_t =
  \ga \sigma^2 \int_0^t \<\phi^2, {\mu}_s \> ds
  \qquad \hbox{for every } t>0 \hbox{ and } \phi \in K_a(\R^d).
 \eeqnn

 Here the filtration of choice is
 ${\cal F}_t:=\sigma\{ \<\phi,\mu_{s}\> , M^\mu_s(\phi), X_s(\phi): \phi \in K_a(\R^d), s \leq t\}$. 
 For any $a\ge0$ and any initial value ${\mu}_{0}\in M_a(\R^d)$, 
this unique law on the Borel subsets of $C([0, \infty), M_a(\R^d))$ 
will henceforth be denoted by $\P_{\mu_0}$ and  
the corresponding expectation by $\E_{\mu_0}$. 
Parameter $\gamma>0$ is related to the branching rate of the particle system
 and $\sigma^2 > 0$ is the variance of the limiting offspring distribution.
 
Almost sure statements, including those referring to (\ref{LT}),  
will henceforth be meant to hold $\P_{\mu_0}$-almost surely.  
 
The properties of law $\P_{\mu_0}$ are determined by way of a function-valued dual process 
for SDSM due to Dawson et al. \cite{DVW19}, a version of the 
original first built by Dawson and Hochberg \cite{DawsonHochberg79} and later generalized 
by Dawson and Kurtz \cite{DawsonKurtz82} as well as others thereafter. 

Let $\{J_t: t\ge 0\}$ be a decreasing c\`{a}dl\`{a}g Markov jump process
on the nonnegative integers $\{0,1,2,\ldots\}$, started at $J_0=m$ 
and decreasing by 1 at a time, 
with Poisson waiting times of intensity $\gamma \sigma^2 l(l-1)/2$
when the process has reached value $l\ge2$. The process is frozen in place  
when it reaches value $1$ and never moves if it is started at 
either $m=0$ or $1$.
Write $\{\tau_k: 0\le k\le J_0-1\}$ for the sequence of jump times of 
$\{J_t:t\ge 0\}$ with $\tau_0=0$ and $\tau_{J_0} = \infty$. 
At each such jump time a randomly chosen projection is effected on the 
function-valued process of interest, as follows. 
Let $\{S_k: 1\le k\le J_0\}$ be a sequence of random
operators which are conditionally independent given $\{J_t: t\ge 0\}$ and satisfy
$$
\P\{S_k = \Phi_{ij}^m | J_{\tau_k -} =m\} = \frac{1}{m(m-1)}, \qquad
1 \le i \neq j \le m,
$$
as long as $m\ge2$. Here $\Phi_{ij}^mf$ is a mapping from $B((\R^d)^{m})$ into $B((\R^d)^{m-1})$ defined by
\beqlb \label{restriction}
\Phi_{ij}^m f(y):= f(y_{1}, \cdots, y_{j-1},y_{i},y_{j+1},\cdots, y_{m}),
\eeqlb
for any $m\ge2$ and $y=(y_1, \cdots, y_{j-1}, y_{j+1}, \cdots, y_m)\in (\R^d)^{m-1}$. 

For some integer $m\ge0$ are given starting values $J_0=m$ 
and $Y_0\in B((\R^d)^{m})$, a bounded function. 
Define process $Y:=\{Y_t: t\ge0\}$, started at $Y_0$ within 
the (disjoint) topological union $\bB:=\cup_{m=0}^\infty B((\R^d)^{m})$, by
\beqlb \label{Yprocess}
Y_t = P^{J_{\tau_k}}_{t-\tau_k} S_k P^{J_{\tau_{k-1}}} _{\tau_k
-\tau_{k-1}}S_{k-1} \cdots P^{J_{\tau_1}}_{\tau_2 -\tau_1} S_1
P^{J_0}_{\tau_1}Y_0, \quad \tau_k \le t < \tau_{k+1}, 0\le k\le J_0-1.
\eeqlb
The process $Y$ is a well-defined $\bB$-valued 
strong Markov process for any starting point $Y_0\in \bB$: it is well-defined on $\bB$  
since all linear operators $\Phi_{ij}^m$ and $P^{m}_{t}$ are contractions on bounded functions; 
and the strong Markov property holds since $Y$ is a pure jump process with finitely many jumps separated 
by exponential waiting times. 
(When $m=0$ we simply write $B((\R^d)^{0})=\R^d$ and $Y_t=P_t^0$ acts as
the identity mapping on constant functions.) 
Clearly, $\{(J_t, Y_t): t\ge 0\}$ is also a strong Markov process. \\

We show next that $Y_t\in\bB$ still holds for all $t>0$ when 
$Y_0\not\in\bB$, under some mild conditions. 
 The following extends several of the results from Dawson et al. \cite{DVW19}. 
 For the proof of the new features, see Subsection \ref{app:pf_th_tha}.

\bth \label{tha}
  Assume that Hypotheses $\ref{hyp:basicassumpFilter}$ 
  and $\ref{hyp:basicassumpElliptic}$ are satisfied. For any $a\ge0$, $m\ge1$,
 $\mu_0\in M_a(\R^d)$, $f\in L^1(\mu_0^m)$ and $t\in[0,\infty)$, there holds
 \beqlb \label{dual}
 \E_{\mu_0} \<f, \mu_t^m\> = \E \left[\<Y_t, \mu_0^{J_t}\>
\exp\left(\frac{\gamma \sigma^2}{2}\ \int_0^t  J_s(J_s-1)ds\right)\Bigl|(J_0,Y_0)=(m,f) \right] 
\eeqlb
with both sides finite in every case. Moreover, for any $p\ge1$, 
every $f \in L^p(\mu_0^m)$ also belongs $\P_{\mu_0}$-almost surely to $L^p(\mu_t^m)$ 
as well as $Y_t \in L^p(\mu_0^{J_t})\cap L^p(\lambda_0^{J_t})\cap C_b^2((\R^d)^{J_t})$, for all $t > 0$.
In addition, for any $T > 0$, there is a constant $c=c(a,d,m,q^m,T) > 0$, independent of the choice of $f$, 
such that  
 \beqlb \label{uniformbound}
 \sup_{0 < t \leq T} \| {\cal I}_{a, J_t}^{-1} Y_t \|_\infty \le c\cdot \| {\cal I}_{a, m}^{-1} f \|_\infty
\eeqlb
and therefore also
 \beqlb \label{Lpmubound}
 \sup_{0 < t \leq T} \<|Y_t|^p, \mu_0^{J_t}\> \le c\cdot \| {\cal I}_{a, m}^{-1} |f|^p \|_\infty
 \cdot \max \left(1, \<I_a,\mu_0\>^m\right)
\eeqlb
both hold $\P_{\mu_0}$-almost surely, whenever the respective right hand side is finite.
\eth
 
 Equation (\ref{dual}) is called the duality identity  
 between the law of SDSM and that of its dual process $(J,Y)$. 


\subsection{Tanaka representation and main result}

For the single particle transition density $q_t^{1}(0,x)$ (from 0) exhibited in 
(\ref{eqn:Semigroup}) for the semigroup $P_t^1$
associated with generator $G_1$ from (\ref{eqn:Gn}), its Laplace transform (in the time variable) is given by  
  \beqlb\label{eqn:green}
 Q^{\lambda}(x) := \int_0^{\infty}e^{-\lambda t} q_t^{1}(0,x)dt,
 \eeqlb
 for any $\lambda > 0$. Formally $Q^{0}$ is known as Green's function for density $q_t^{1}$ 
 and exhibits a potential singularity at $x=0$. 
 By (\ref{LSU}), for all $x \in \R^d \smallsetminus \{0\}$ we can also write
   \beqlb\label{eqn:interchange}
    \partial_{x_i}Q^{\lambda}(x)= 
    \partial_{x_i}\int_0^{\infty}e^{-\lambda t}q_t^{1}(0,x)dt
    =\int_0^{\infty}e^{-\lambda t}\partial_{x_i}q_t^{1}(0,x)dt < \infty
 \eeqlb
  for any $i \in \{1,2,\cdots,d\}$, with the derivative taken in the distributional sense. 
  
When $d=1$, under Hypothesis \ref{hyp:basicassumpElliptic}, 
there holds $ \| I_{a}^{-1} (Q^{\lambda})^p \|_\infty < \infty$ for any $a\ge0$ and $p\ge1$ 
--- see inequality (\ref{gamma}) in Subsection \ref{app:pf_prop_mu0integrability} --- 
so the bounds (\ref{uniformbound}) and (\ref{Lpmubound}) are applicable in that case 
(though it is not the one of interest here). Unfortunately, when $d\ge2$, 
even in the special case of the gaussian kernel $\varphi$, there holds 
$ \| Q^{\lambda} \|_\infty = \infty$. Nevertheless, the $\|\cdot\|_{\mu_0,p}$-norms,  
with $p\in[1,\infty)$ instead, yield preliminary bounds. 
 
 \bprop \lab{prop_mu0integrability}
 Under Hypotheses $\ref{hyp:basicassumpFilter}$   
 and $\ref{hyp:basicassumpElliptic}$, for every $\lambda > 0$ and every measure 
 $\mu_0\in\cup_{a\ge0}M_a(\R^d)$ satisfying Hypotheses $\ref{hyp:basicassumpGaussUniform}$ 
 and $\ref{hyp:basicassumpUniformInteg}$, 
 with constants $A^*>0$ and $c>0$ from $(\ref{Aronsonbounds})$, 
 \[ 
 \sup_{0 < s \leq t}\sup_{w\in \R^d} \E_{\mu_0} \<(Q^{\lambda}(w-\cdot))^p, \mu_s\> 
     \leq A^*\Upsilon_{ct}(\varphi, \mu_0)\<(Q^{\lambda})^p, \lambda_0\> < \infty 
\]  
 holds at any $t>0$ with 
\[ 
\sup_{w\in \R^d} \<(Q^{\lambda}(w-\cdot))^p, \mu_0\> < \infty 
\] 
at $t=0$ and so does $Q^{\lambda} \in L^p(\mu_0)\cap L^p(\lambda_0)$, 
in each of the following three cases: \\
  {\em (i)} when $p=1$, for all $d\ge1$; \\ 
  {\em (ii)} when $p=2$ and $d\le3$; \\
  {\em (iii)} for all $p\ge1$ when $d=1$ or $2$. \\ 
Furthermore, there holds \\ 
  {\em (iv)} in all dimensions $d\ge1$ and for any $i \in \{1,2,\cdots,d\}$, 
\[ 
\sup_{0 < s \leq t}\sup_{w\in \R^d} \E_{\mu_0} \int_{\R^d}|\partial_{x_i}Q^{\lambda}(w-x)| \mu_s(dx) 
     \leq A^*\Upsilon_{ct}(\varphi, \mu_0)\<|\partial_{x_i}Q^{\lambda}|, \lambda_0\> < \infty 
\] 
at any $t>0$ with 
\[ 
\sup_{w\in \R^d} \int_{\R^d}|\partial_{x_i}Q^{\lambda}(w-x)| \mu_0(dx)< \infty 
\] 
at $t=0$ and so does $\partial_{i}Q^{\lambda} \in L^1(\mu_0)\cap L^1(\lambda_0)$.
\eprop

In the special case of Lebesgue measure $\mu_0=\lambda_0\in M_a(\R^d)$ (when $a > d$), 
the pointwise part of each of statements (i), (ii) and (iv) pertaining to $t=0$ was obtained in \cite{DVW19}, where  
$\partial_{x}Q^{\lambda} \in L^2(\lambda_0)$ is also proved to hold when $d=1$ and to fail when $d\ge2$. 
This pointwise finiteness suffices to prove the existence of the local time under Hypotheses 
\ref{hyp:basicassumpFilter}, \ref{hyp:basicassumpElliptic}, \ref{hyp:basicassumpGaussUniform} and  
\ref{hyp:basicassumpUniformInteg}; but the uniform extensions presented here are useful in 
our proof of its joint H{\"o}lder continuity. The proof of Proposition \ref{prop_mu0integrability} 
is found in Subsection \ref{app:pf_prop_mu0integrability}.  
Our main result can now be stated. 

  \bth \lab{lt_th1} 
 Under Hypotheses $\ref{hyp:basicassumpFilter}$   
 and $\ref{hyp:basicassumpElliptic}$, with $d=1$, $2$ or $3$, 
 select any $a\ge0$ and ${\mu}_{0}\in M_a(\R^d)$ satisfying both  
 Hypotheses $\ref{hyp:basicassumpGaussUniform}$ and $\ref{hyp:basicassumpUniformInteg}$, 
 with ambiant Brownian sheet $W$ 
 and martingale measure $M^\mu$. 
  For every choice of $(t,x) \in [0, \infty) \times \R^d$, 
   \beqlb \lab{TanakaI}
    \Lambda^x_{t} & := & \<Q^{\lambda}(x - \cdot), \mu_0\> - \<Q^{\lambda}(x- \cdot),
    \mu_{t}\>
    + \lambda \int_0^{t} \<Q^{\lambda}(x - \cdot), \mu_s\>ds \nonumber \\
    & &  + \sum_{p=1}^{d}\int_0^{t} \int_{\R^d}\<h_p(y- \cdot)
    \partial_p Q^{\lambda}(x - \cdot), \mu_s\>W(dy,ds) \nonumber \\
   & &  + \int_0^{t} \int_{\R^d}Q^{\lambda}(x-y)M^\mu(dy,ds)
   \eeqlb
  is the local time for SDSM $\{\mu_t\}$.  
  It satisfies both $(\ref{LT})$ $\P_{\mu_0}$-almost surely for every choice of $\phi \in C_c(\R^d)$ 
  and square-integrability for every $(t,x) \in [0, \infty) \times \R^d$ --- in fact there holds  
  \beqlb \lab{Tanaka2}
 \sup_{x\in\R^d}\sup_{0 \leq t \leq T} \E_{\mu_0} \left[ \left| \Lambda^x_{t} \right|^2\right] < \infty.
 \eeqlb 
 Moreover, there exists a version of $\Lambda^x_{t}$ 
 which is H\"{o}lder jointly continuous in $(t,x) \in [0, \infty) \times \R^d$. Explicitly, 
 for every exponent $\alpha\in(0,1)$ in space, as well as for every choice of $n > 1$, 
 there is a constant $c=c(T, d, h, q^1, n, \alpha, \lambda)> 0$ such that,
for every $x,z\in\R^d$ and $0 < s < t \le T$ there holds
 \beqlb \lab{JointHolder}
 \left[ \E_{\mu_0} \left| \Lambda^z_{t}-\Lambda^x_{s}\right|^{2n} \right]^{1/2n} 
 \leq c\cdot (|z - x|^{\alpha}+|t - s|^{1/2}). 
 \eeqlb 
 \eth

This last inequality suffices for an application of Kolmogorov's continuity criterion 
(such as Theorem I.2.1 of Revuz and Yor \cite{RY94}), which  
for each compact $K \subset [0, \infty) \times \R^d$, states the existence of a 
modification $\tilde\Lambda^x_{s}$ of $\Lambda^x_{s}$ such that 
\[
\E_{\mu_0} \left[ \sup_{(s,x),(t,z)\in K \atop (s,x)\neq(t,z)}\left\{
\frac{\left| \tilde\Lambda^z_{t}-\tilde\Lambda^x_{s}\right|}{(|z - x|^{\gamma}+|t - s|^{\beta})} 
\right\}^{2n}\right] <\infty
\] 
holds for every $\gamma\in[0,1)$ and $\beta\in[0,1/2)$. 

Equation (\ref{TanakaI}) is from Dawson et al. \cite{DVW19} 
and is called the Tanaka formula for SDSM.  Its right 
hand side is a Schwartz distribution and (\ref{LT}) holds in the distributional sense.  
The value of the local time is independent of $\lambda > 0$ but it varies with $d$. 
What is new here of course is the last assertion, regarding joint H{\"o}lder continuity. 
Its proof occupies the next section, where we also explore the impact of 
Hypothesis $\ref{hyp:basicassumpUniformInteg}$. 


  \section{Joint continuity of the local time}\label{sec:Tanaka} 
  \setcounter{equation}{0}

 The proof of Theorem \ref{lt_th1} proceeds next, using some of the sharp inequalities 
 previously mentionned, by comparing the higher moments of both the local time for 
 SDSM $\{\mu_t\}$ and its space-time differences, against those associated with a 
 special function controlling the growth of the two martingales emanating from the branching 
 and random environment mechanisms within the Tanaka formula (\ref{TanakaI}). 
 The argument is in the spirit of a similar result by 
 Ethier and Krone \cite{EthierKrone95}) for the Fleming-Viot process, for which they use 
 Super-Brownian motion as the benchmark. 

Denote by $\tilde{Q}^{\lambda}$ the Laplace transform in the special case 
 of Super-Brownian motion, where $h\equiv0$ and $c$ is the identity matrix: 
 \beqlb\label{eqn:greenSBM}
 \tilde{Q}^{\lambda}(x) := \int_0^{\infty}e^{-\lambda s} \varphi_s(x)ds .
 \eeqlb

 In view of the components in the Tanaka formula (\ref{TanakaI}), 
 the key functions to bound are $(\tilde{Q}^{\lambda}(x))^2$ and the following, 
 defined for every $\lambda>0$:
 \beqlb \lab{GaussianKey}
 \tilde{G}^{\lambda}(x) := \int_0^{\infty} e^{-\lambda s} s^{-1/2} \varphi_s(x) ds .
 \eeqlb

The two are linked by way of the following technical result. 

 \ble \lab{le_SquareVSPartial}
 Fix $\lambda>0$ and $d\ge1$. The following hold for all $x\in\R^d$: \\
  {\em (i)} $\tilde{Q}^{\lambda}(x)\le \max(1,\lambda^{-1/2})\cdot\tilde{G}^{\lambda/2}(x)$;  \\
  {\em (ii)} $(\tilde{Q}^{\lambda}(x))^2 \le K \tilde{G}^{2\lambda}(x)$ for some constant 
  $K=K(d,\lambda)>0$, provided $d\le3$.
  \ele
\proof 
Elementary inequality (\ref{gamma}) applied to   
$\tilde{Q}^{\lambda}(x)\le \tilde{G}^{\lambda/2}(x)\sup_{s\ge0} e^{-\lambda s/2} s^{1/2}$ yields the first statement. 
For every $u>0$ and $v>0$, there holds 
 \beqlb \lab{gaussianproduct}
\varphi_u(x)\varphi_v(x) = (2\pi w)^{- d/2} \varphi_{z}(x)
 \eeqlb 
 with $z=uv/(u+v)$ and $w=u+v$, as seen 
by expanding the right side. 
The Jacobian $J(w,z;u,v)=|v-u|/(u+v)=\sqrt{1-4z/w}$ of this change of variable 
is non null except on the line $v=u$.  
Keeping in mind that, for all $w>0$ and $z>0$, each dimension amongst $d=1$, $2$ or $3$ 
satisfies $w^{-d/2}\le w^{-3/2}+w^{-1/2}$ and $z^{1-d/2}\le z^{-1/2}+z^{1/2}$ everywhere, 
there is a constant $K=K(d,\lambda)>0$ such that 
 \beqlb \lab{Jacobian}
 (\tilde{Q}^{\lambda}(x))^2 & = & 2\int_0^{\infty}dz \int_{4z}^{\infty} 
 e^{-\lambda w}  (2\pi w)^{-d/2} \varphi_z(x) J(w,z;u,v)^{-1}dw  \nonumber \\
  & = & 2\int_0^{\infty}\varphi_z(x)dz \int_{4z}^{\infty} 
  e^{-\lambda w}(2\pi w)^{-d/2}(1-4z/w)^{-1/2}dw  \nonumber \\
 & \le & 2(8\pi)^{-d/2}\int_0^{\infty}e^{-4\lambda z}z^{-d/2}\varphi_z(x)dz 
 \int_{4z}^{8z}(1-4z/w)^{-1/2}dw  \nonumber \\
& + & 2^{3/2}(2\pi)^{-d/2}\int_0^{\infty}\varphi_z(x)dz \int_{8z}^{\infty} 
  e^{-\lambda w}w^{-d/2}dw  \nonumber \\
& \le & 8(8\pi)^{-d/2}\int_0^{\infty}e^{-4\lambda z}z^{1-d/2}\varphi_z(x)dz \int_{0}^{1/2} 
s^{-1/2}(1-s)^{-2}ds  \nonumber \\
& + & 2^{3/2}(2\pi)^{-d/2}\int_0^{\infty}\varphi_z(x)dz \int_{8z}^{\infty} 
  e^{-\lambda w}(w^{-3/2}+w^{-1/2})dw  \nonumber \\
& \le & 64(8\pi)^{-d/2}\int_0^{\infty}e^{-4\lambda z}z^{1-d/2}\varphi_z(x)dz  \nonumber \\
& + & (2+\lambda^{-1})(2\pi)^{-d/2}\int_0^{\infty}e^{-8\lambda z}z^{-1/2}\varphi_z(x)dz \nonumber \\ 
& \le & K \left( \tilde{G}^{4\lambda}(x) + \tilde{G}^{2\lambda}(x) + \tilde{G}^{8\lambda}(x) \right)
 \eeqlb 
 after splitting the second integral, taking advantage of the monotonicity of the various functions, 
 using the change of variable $s=1-4z/w$ with $dw=4zds/(1-s)^2$ and integrating by parts. 
 Finally note that inequality (\ref{gamma}) yields 
\[
\int_0^{\infty}e^{-4\lambda z}z^{1/2}\varphi_z(x)dz 
\le  \tilde{G}^{2\lambda}(x) \cdot\sup_{z\ge0}ze^{-2\lambda z} 
 \le  \tilde{G}^{2\lambda}(x) \cdot\max(1,1/2\lambda)
\] 
and that $\tilde{G}^{\lambda}(x)$ is decreasing in $\lambda$ at every $x\in\R^d$. 
 \qed
 
 \noindent  {\bf Remark:} 
 Under Hypothesis \ref{hyp:basicassumpElliptic}, 
 Aronson's inequalities (\ref{Aronsonbounds}) ensure that both statements remain valid 
 when the specific $\tilde{Q}^{\lambda}(x)$ is replaced by the general $Q^{\lambda}(x)$, 
 provided the constants in front of $\tilde{G}$, as well as the scale parameters inside it, 
 are adjusted accordingly. The same is true when bounding either $|\partial_{x_i}Q^{\lambda}(x)|$ 
 or $|\partial_{x_i}\tilde{Q}^{\lambda}(x)|$ above by a scaled multiple of $\tilde{G}$, 
 using inequality (\ref{LSU}) instead. The key inequality in Lemma \ref{le_SquareVSPartial} is of 
 course the second one. It controls the size and fluctuations of the square of 
$\tilde{Q}^{\lambda}(x)$, which drives the quadratic variation process for the branching martingale, 
by way of the more volatile size and fluctuations of $\tilde{G}^{\lambda}(x)$, which drives the 
quadratic variation process for the random environment martingale. 

  \ble \lab{le_mutintegrability}
 Assume that Hypotheses $\ref{hyp:basicassumpFilter}$, $\ref{hyp:basicassumpElliptic}$,  
 $\ref{hyp:basicassumpGaussUniform}$ and $\ref{hyp:basicassumpUniformInteg}$ 
 are satisfied for some $a\ge0$ and $\mu_0\in M_a(\R^d)$. 
 For any $t>0$, $b>0$, $\lambda > 0$, $d\le3$ and $n\ge1$, there holds 
 \beqlb \lab{mutintegrabilityI}
 \sup_{0 < u \leq t}\sup_{x\in \R^d}  \E_{\mu_0} \<I_a(x-\cdot), \mu_u\>^{n} < \infty ;
 \eeqlb
\beqlb \lab{mutintegrabilityU}
  \sup_{0 < u \leq t}\sup_{x\in \R^d}  \E_{\mu_0}
  \left\{\left[\int_0^{b} s^{-1/2}\<\varphi_{s}(x-\cdot), \mu_u\>ds\right]^{n}\right\}  < \infty ;
 \eeqlb
\beqlb \lab{mutintegrabilityG}
 \sup_{0 < u \leq t}\sup_{x\in \R^d}  \E_{\mu_0} \<\tilde{G}^{\lambda}(x-\cdot), \mu_u\>^{n} < \infty ;
 \eeqlb
 \beqlb \lab{mutintegrabilityQ}
\sup_{0 < u \leq t}\sup_{x\in \R^d}  \E_{\mu_0} \< [\tilde{Q}^{\lambda}(x-\cdot)]^2, \mu_u\>^{n} < \infty .
 \eeqlb
\ele
  \proof 
  We first prove (\ref{mutintegrabilityI}). 
For any Borel measurable, real-valued function $f$ on some $(\R^d)^n$ with $n\ge1$, 
define mapping $T_X^nf(Y)=f(X-Y)$ and observe the following three properties: 
$T_{X_{n-1}}^{n-1}\circ\Phi_{ij}^n=\Phi_{ij}^n\circ T_{X_n}^n$ holds everywhere with $X_n=(x,\ldots,x)\in(\R^d)^{n}$; 
while $T_{X_n}^n$ does not necessarily commute with $P_u^n$ in general, it does so when the kernel $q^1$ is 
symmetric, in particular when $q_u^{1}(\cdot,y)=\varphi_u(y-\cdot)$, so we can write 
$T_{X_n}^n\circ \tilde{P}_u^n=\tilde{P}_u^n\circ T_{X_n}^n$ with $\tilde{P}^n$ the heat semigroup on $(\R^d)^n$; 
and finally, writing the dual as $Y_u=Y_u(f)$ to indicate the dependency on $Y_0=f$ 
and $\tilde{Y}_u=\tilde{Y}_u(f)$ for its special case in (\ref{Yprocess}) involving $\tilde{P}_{cu}^n$ instead of $P_u^n$, 
that is, the dual to a rescaled Super-Brownian motion, upper bound (\ref{Aronsonbounds}) 
applied to $q^n$ instead of $q^1$ implies 
\[ 
Y_u\circ T_{X_n}^n(f)\le A^*\cdot\tilde{Y}_{u}\circ T_{X_n}^n(f)=A^*\cdot T_{X_{J_u}}^{J_u}\circ \tilde{Y}_{u}(f).
\] 
Put $C_1=\exp(\gamma \sigma^2nt(n-1)/2)$ from (\ref{dual}). 
Using bounds (\ref{uniformbound}) and (\ref{Lpmubound})  
with $f={\cal I}_{a, n}$ and constant $C_2$ to differentiate it from $c$ already in use here, we get 
  \[ 
 \E_{\mu_0} \<I_a(x-\cdot), \mu_u\>^{n} 
 \le C_1A^*\cdot\E\<\tilde{Y}_{u}({\cal I}_{a, n}), \mu_0(x-\cdot)^{J_u}\>
 \le C_1C_2A^*\max \left(1, \<I_a,\mu_0(x-\cdot)\>^n\right),
\] 
uniformly in both $x\in\R^d$ and $u\in(0,t]$, 
and hence (\ref{mutintegrabilityI}) follows from Hypothesis \ref{hyp:basicassumpUniformInteg}.

To prove (\ref{mutintegrabilityU}), H\"{o}lder's inequality yields, for any $u\in(0,t]$, $b>0$ and $n\ge2$, 
 \[ 
\left[\int_0^{b} s^{-1/2}\<\varphi_{s}(x-\cdot), \mu_u\>ds\right]^{n} 
\le \left[\int_0^b s^{\alpha} \<\varphi_{s}(x-\cdot), \mu_u\>^{n}ds\right]
            \left[\int_0^b \left(s^{-1/2}s^{-\alpha/n}\right)^{[n/(n-1)]}ds \right]^{n-1} 
\] 
with the last integral finite as soon as $\alpha<(n-2)/2$. Therefore, to get (\ref{mutintegrabilityU}) 
it suffices to prove the existence of a small $\epsilon>0$ such that the sequence $\alpha_n=-\epsilon+(n-2)/2>-1$ 
for $n\ge2$, starting at $\alpha_1=-1/2$, verifies, for all $n\ge1$, 
\beqlb \lab{mutintegrabilityU2}
  \sup_{0 < u \leq t}\sup_{x\in \R^d}  \E_{\mu_0}
  \left\{\int_0^{b} s^{\alpha_n}\<\varphi_{s}(x-\cdot), \mu_u\>^{n}ds\right\}  < \infty .
 \eeqlb 
Just as in the proof of (\ref{mutintegrabilityI}), but this time under Hypothesis \ref{hyp:basicassumpGaussUniform}, write 
 \[ 
 \E_{\mu_0} \left[\<\varphi_{s}(x-\cdot), \mu_u\>^{n}\right] 
 = \E_{\mu_0} \<f_n(X_n-\cdot), \mu_u^{n}\> 
 \le C_1A^*\cdot\E\<\tilde{Y}_{u}(f_n), \mu_0(x-\cdot)^{J_u}\> 
\] 
where we set $f_n(x_1,\ldots,x_n)=\prod_{k=1}^n\varphi_{s}(x_k)$ with $0<s\le b$. 

For any $n\ge1$, in the event $\{\tau_1>u\}$ where the first jump does not occur prior to time $u>0$, 
 $\tilde{Y}_{u}(f_n) =\tilde{P}^n_u f_n(x_1,\ldots,x_n)
 =\prod_{k=1}^n\tilde{P}^1_u \varphi_{s}(x_k)=\prod_{k=1}^n\varphi_{u+s}(x_k)$ 
holds by the Markov property, so there ensues, for any $b>0$ and $0<s\le b$, 
 \[  
 \sup_{0 < u \leq t}\sup_{x\in \R^d} \E \left\{\<\tilde{Y}_{u}(f_n), \mu_0(x-\cdot)^{J_u}\>\bigg| \tau_1>u \right\}  
\le \left[\Upsilon_{c(t+b)}(\varphi, \mu_0)\right]^{n}
 \]
and hence, keeping in mind that $f_n$ depends on $s$,  
 \[  
  \sup_{0 < u \leq t}\sup_{x\in \R^d} \E 
  \left\{\left[\int_0^{b} s^{\alpha_n}\<\tilde{Y}_{u}(f_n), \mu_0(x-\cdot)^{J_u}\>ds\right] \bigg| \tau_1>u \right\} 
  \le \left[\Upsilon_{c(t+b)}(\varphi, \mu_0)\right]^{n}
     \int_0^b s^{\alpha_n}ds,
 \]
which is finite by Hypothesis \ref{hyp:basicassumpGaussUniform}. 
This also completes the proof of (\ref{mutintegrabilityU2}) in the case $n=1$.

Next, with $n\ge2$, immediately after the first jump time $\tau_1$ of the function-valued dual,  
when some randomly selected coordinate $j$ is replaced by some other $i$, but before the second jump, 
hence in $\{\tau_1 < u < \tau_2\}$ and for any $\epsilon>0$, we get 
 \beqlb \lab{onejump1}
& & \tilde{P}^{n-1}_\epsilon \tilde{Y}_{\tau_1}(f_n)(x_1,\ldots,x_{j-1},x_{j+1},\ldots,x_n)  \nonumber \\ 
& = & \biggl[\prod_{k=1\atop k\neq i,j}^n\varphi_{\epsilon+\tau_1+s}(x_k)\biggr] 
\tilde{P}^1_\epsilon(\varphi_{\tau_1+s}^2)(x_i) \nonumber \\ 
& = & \biggl[\prod_{k=1\atop k\neq i,j}^n\varphi_{\epsilon+\tau_1+s}(x_k)\biggr] 
 [4\pi (\tau_1+s)]^{- d/2} \varphi_{\epsilon+(\tau_1+s)/2}(x_i)
 \eeqlb 
where we used $\varphi_{u}^2(x_i)= (4\pi u)^{- d/2} \varphi_{u/2}(x_i)$ from (\ref{gaussianproduct}) with $v=u$. 
Trajectory by trajectory, the function-valued dual is therefore bounded everywhere in $\{\tau_1 < u < \tau_2\}$ by 
 \[ 
 \<\tilde{Y}_{u}(f_n), \mu_0(x-\cdot)^{J_u}\> 
\le [4\pi (\tau_1+s)]^{-d/2} \left\< \varphi_{u-\tau_1+(\tau_1+s)/2}, \mu_0(x-\cdot)\right\>
\left[\Upsilon_{c(t+b)}(\varphi, \mu_0)\right]^{J_u-1}
 \] 
for the selected pair $(i,j)$ with $J_u=n-1$ and $n\ge2$. Choosing amongst ${n\choose2}$ such pairs  
with equiprobability, there holds almost surely
 \beqlb \lab{onejump2}
& &   \sup_{0 < u \leq t}\sup_{x\in \R^d}  \E 
\left\{\left[\int_0^{b} s^{\alpha_n} \<\tilde{Y}_{u}(f_n), \mu_0(x-\cdot)^{J_u}\> ds\right]\bigg| \tau_1 < u < \tau_2\right\} \\ 
  & \le & \left[\Upsilon_{c(t+b)}(\varphi, \mu_0)\right]^{n-1} \int_0^b s^{\alpha_n} [4\pi (\tau_1+s)]^{-d/2}ds 
 < \infty,  \nonumber 
 \eeqlb 
for all $n\ge2$, since $0<\tau_1 < u \le t$ and $s\le b$ yields $u-\tau_1+(\tau_1+s)/2\le t+b$. 
This completes the proof of (\ref{mutintegrabilityU2}) in the case $n=2$ provided $d\le3$ since under that restriction there holds
$\int_0^b s^{\alpha_2} \E\{(\tau_1+s)^{-d/2}\}ds<\infty$ with $\alpha_2=-\epsilon$ for any $\epsilon\in(0,1/2)$. 
(We keep the same $\epsilon$ for the rest of the proof for the sake of simplicity but it is not necessary to do so.)

In $\{\tau_2 < u < \tau_3\}$, either the second jump involves two new coordinates, generating a product of two factors 
$[4\pi (\tau_1+s)][4\pi (\tau_2+s)]$ inside (\ref{onejump1}) instead of just one; 
or coordinate $i$ is involved in both jumps, to or from a new site $k$, for the second jump, 
yielding a different product of the form $[4\pi (\tau_1+s)][2\pi (\tau_2+s+\tau_2-\tau_1+(\tau_1+s)/2)]$, 
by appealing directly to (\ref{gaussianproduct}). Explicitly, with $J_u=n-2$ and $n\ge3$,  
 \beqlb \lab{twojumps1} 
 & &  \<\tilde{Y}_{u}(f_n), \mu_0(x-\cdot)^{J_u}\>  
 \le [4\pi (\tau_1+s)]^{-d/2}  [4\pi (\tau_2+s)]^{-d/2} 
\left[\Upsilon_{c(t+b)}(\varphi, \mu_0)\right]^{J_u} \nonumber \\ 
  & + & 2 [4\pi (\tau_1+s)]^{-d/2}  [2\pi (\tau_2+s+\tau_2-\tau_1+(\tau_1+s)/2)]^{-d/2} 
\left[\Upsilon_{c(t+b)}(\varphi, \mu_0)\right]^{J_u},  \nonumber 
 \eeqlb 
 using the fact that $\min(u,v)/2\le z\le\max(u,v)/2$ holds in (\ref{gaussianproduct}) and ensures that 
 each jump replaces the product $\varphi_u\varphi_v$ with $\max(u,v)\le t+b$ by $\varphi_z$ with $z\le t+b$ every time.
 
Using a generic constant $C>0$ and dropping some of the $\tau_1$ and $\tau_2$ terms in the denominator 
to keep the expression as simple as possible, (\ref{onejump2}) becomes 
 \beqlb \lab{twojumps2}
& &   \sup_{0 < u \leq t}\sup_{x\in \R^d}  \E 
\left\{\left[\int_0^{b} s^{\alpha_n} \<\tilde{Y}_{u}(f_n), \mu_0(x-\cdot)^{J_u}\> ds\right]\bigg| \tau_2 < u < \tau_3\right\} \nonumber \\ 
  & \le & C \left[\Upsilon_{c(t+b)}(\varphi, \mu_0)\right]^{n-2} \int_0^b s^{\alpha_n} 
   (\tau_1+s)^{-d/2} (\tau_2+s)^{-d/2} ds < \infty,  \nonumber 
 \eeqlb 
for all $n\ge3$ and of course $0<\tau_1 < \tau_2 < u \le t$. 
The case $n=3$ of (\ref{mutintegrabilityU2}) ensues since there holds
$\int_0^b s^{\alpha_3} \E\{(\tau_1+s)^{-d/2}(\tau_2+s)^{-d/2}\}ds<\infty$, again under $d\le3$. 

In general, in $\{\tau_m < u < \tau_{m+1}\}$, when $m\ge1$ jumps have occurred, 
the following simple pattern emerges along each trajectory: 
the $n-m$ surviving coordinates amongst the initial $n$, each contribute a factor 
$\Upsilon_{c(t+b)}(\varphi, \mu_0)$ due to the recurring use of (\ref{gaussianproduct}), which ensures that 
each jump replaces the product $\varphi_u\varphi_v$ with $\max(u,v)\le t+b$ by $\varphi_z$ with $z\le t+b$; 
the first jump contributes a factor of $(\tau_1+s)^{-d/2}$, the second one a factor smaller than $(\tau_2+s)^{-d/2}$, 
because of the resulting sum $w=u+v$ in (\ref{gaussianproduct}), the third one a factor smaller than $(\tau_3+s)^{-d/2}$, 
for the same reason, and so on. As there are finitely many choices of coordinates from the dual mechanism, we can write, 
with another generic constant $C>0$, 
 \beqlb \lab{mjumps2}
& &   \sup_{0 < u \leq t}\sup_{x\in \R^d}  \E 
\left\{\left[\int_0^{b} s^{\alpha_n} \<\tilde{Y}_{u}(f_n), \mu_0(x-\cdot)^{J_u}\> ds\right]\bigg| \tau_m < u < \tau_{m+1}\right\} \nonumber \\ 
  & \le & C \left[\Upsilon_{c(t+b)}(\varphi, \mu_0)\right]^{n-m} \int_0^b s^{\alpha_n} ds 
  \prod_{i=1}^{m} (\tau_i+s)^{-d/2} < \infty,  \nonumber 
 \eeqlb 
for all $n\ge m+1$ and $0<\tau_1 < \tau_2 < \dots < \tau_m < u \le t$. 
The case $n=m+1$ of (\ref{mutintegrabilityU2}) ensues since there holds
\[
\int_0^b s^{\alpha_{m+1}} \E\left\{ \prod_{i=1}^{m} (\tau_i+s)^{-d/2} \right\}ds<\infty,
\]
again under $d\le3$ and the proofs of (\ref{mutintegrabilityU}) and (\ref{mutintegrabilityU2}) are complete. 

That (\ref{mutintegrabilityI}) and (\ref{mutintegrabilityU}) together imply (\ref{mutintegrabilityG}) 
is a consequence of inequalities relegated to Subsection \ref{app:pf_prop_mu0integrability} 
and already used in the proof of Proposition \ref{prop_mu0integrability}, where starting with 
(\ref{pnorm2}) we used the supremum $\Upsilon$ from (\ref{nonBAUnif}) in order to handle 
all at once the powers $p\ge1$ required in the proof of Proposition \ref{prop_mu0integrability}.
Going back to (\ref{pnorm2}) with $p=1$ yields a sharper version of (\ref{partialboundat0}), 
which nevertheless remains valid when substituting $\tilde{G}$ for $|\partial_{x_i}Q^{\lambda}(x)|$ 
and setting values $a_1=1$, $a_2=1/2$, $\eta=2a_2p/a=1/a$ and $\eta^*=\eta/2a_2=\eta=1/a$. 
We obtain instead
\[ 
\<\tilde{G}^{\lambda}(x-\cdot), \mu_u\> \le C_1 \<I_a(x-\cdot), \mu_u\>
+ a_1(2\pi)^{d/2} \int_0^\eta s^{-1/2} \<\varphi_{s}(x-\cdot), \mu_u\> ds
\] 
for all $x\in\R^d$ and using any $a>0$ in the base case where $\mu_0\in M_0(\R^d)$. 

Finally, by Lemma \ref{le_SquareVSPartial}, (\ref{mutintegrabilityG}) implies (\ref{mutintegrabilityQ}). 
  \qed 
   
  \bth \lab{th_HolderEnvironment} 
  Assume that Hypotheses $\ref{hyp:basicassumpFilter}$, $\ref{hyp:basicassumpElliptic}$,  
 $\ref{hyp:basicassumpGaussUniform}$ and $\ref{hyp:basicassumpUniformInteg}$ 
  are satisfied for some $a\ge0$ and $\mu_0\in M_a(\R^d)$. For either $d=1$, $2$ or $3$, the random field 
 \beqlb \lab{DI0}
  \Xi_t(x):= \int_0^{t} \int_{\R^d} \<h_p(y-\cdot) \partial_pQ^{\lambda}(x-\cdot),
 \mu_s\> W(dy,ds) 
 \eeqlb
is a square-integrable ${\cal F}_t$-martingale, 
for every $\lambda > 0$ and $p\in\{1,2,\ldots,d\}$, with quadratic variation given by
\beqnn
  \<\Xi(x)\>_t 
  = \int_0^{t} ds \int_{\R^d} \<h_p(y-\cdot)\partial_pQ^{\lambda}(x-\cdot), \mu_s\>^2 dy  
 \eeqnn
 and satisfying $\sup_{x\in\R^d} \E_{\mu_0}  \<\Xi(x)\>_t  < \infty $ for every $t > 0$.
 There exists a version of $\{\Xi_t(x), t \geq 0, x \in \R^d \}$ which is 
jointly H\"{o}lder continuous in $(t,x) \in [0, \infty) \times \R^d $,
for every $\lambda > 0$ and $p\in\{1,2,\ldots,d\}$.
  \eth

  \proof 
 The statement regarding (\ref{DI0}) is from Dawson et al. \cite{DVW19} so we focus on the joint H\"{o}lder continuity. 
 Given any suitably measurable real-valued functions $H\ge0$ and $F\ge0$, there holds,
 for any $z\in\R^d$, $n\ge1$ and $0<s<t<\infty$,
 \beqlb \lab{CauchySchwarzHolder}
 & & \left[ \int_{s}^{t}du \int_{\R^d}\lambda_0(dy) \left\{\int_{\R^d} H(y,w) F(x,z,w) \mu_u(dw) \right\}^2\right]^n
 \nonumber \\
 \le & & \left(\sup_{w_1,w_2}\int_{\R^d}H(y,w_1)H(y,w_2)\lambda_0(dy)\right)^n 
            \left(\int_{s}^{t} \<F(x,z,\cdot), \mu_u\>^{2n}du\right) 
            \left(\int_{s}^{t}du\right)^{n-1}
 \eeqlb
 by successively expanding the square, applying Fubini's theorem and using H\"{o}lder's inequality 
 with conjugate exponents $n$ and $n/(n-1)$, treating the case $n=1$ separately. 
 
From (\ref{CauchySchwarzHolder}) we get upper bounds for both the time lag
 \beqlb \lab{qvntimediff}
 \<\Xi_\cdot(x)-\Xi_s(x)\>_t^{n} \le {\bar h}^n(t-s)^{n-1}\int_{s}^{t} \< |\partial_pQ^{\lambda}(x-\cdot)|, \mu_u\>^{2n}du  
 \eeqlb
 and the space lag for the quadratic variation process 
  \beqlb \lab{qvnspacediff}
  \<\Xi(z)-\Xi(x)\>_t^n \le {\bar h}^nt^{n-1}\int_0^{t} 
  \< |\partial_pQ^{\lambda}(z-\cdot)-\partial_pQ^{\lambda}(x-\cdot)|, \mu_u\>^{2n}du  
 \eeqlb
after Hypothesis \ref{hyp:basicassumpElliptic} first takes care of the finiteness of constant   
 \beqnn
   {\bar h}:=\max_p\sup_{w_1,w_2}\int_{\R^d} | h_p(y-w_1) | \cdot | h_p(y-w_2) | dy 
   \leq \max_p(\|h_p\|_\infty \cdot \|h_p\|_{\lambda_0,1}) < \infty.
 \eeqnn

Combining (\ref{LSU}) and (\ref{formulaI}), for any $\alpha\in(0,1)$, the constants 
$c=c(\alpha) > 0$ and $c_0=c_0(\alpha) > 0$ from (\ref{formulaI}) can be respectively 
increased and decreased, so that both 
$| \partial_pQ^{\lambda}(x)| \leq c\cdot G^{\lambda^*}(x)$ and 
$| \partial_pQ^{\lambda}(z) -  \partial_pQ^{\lambda}(x)| \leq c|z - x|^{\alpha} 
\left[G^{\lambda^*}(x)+G^{\lambda^*}(z)\right]$ 
 are satisfied, for every $x,z\in {\R}^d$, with $G^{\lambda^*}$ from (\ref{GaussianKey}) rewritten as
 \beqnn
 G^{\lambda^*}(x) = (\frac{1}{2c_0})^{1/2}(\frac{c_0}{\pi})^{d/2}  
          \int_0^{\infty} e^{-\lambda s} s^{- (d+1)/2} \exp{\{-c_0\frac{|x|^2}{s}\}} ds
 \eeqnn
after putting $u=s/2c_0$ and $\lambda^*=2c_0\lambda$ to standardize the exponential kernel 
and get rid of $c_0$, by incorporating it into both $\lambda^*$ and $c$. 
This transforms (\ref{qvntimediff}) and (\ref{qvnspacediff}) into  
 \beqlb \lab{qvntimediff2}
 \<\Xi_\cdot(x)-\Xi_s(x)\>_t^{n} \le c(t-s)^{n-1}\int_{s}^{t} \< G^{\lambda^*}(x-\cdot), \mu_u\>^{2n}du  
 \eeqlb
 with a new $c=c(T, d, h, q^1, n, \alpha)>0$ and  
  \beqlb \lab{qvnspacediff2}
  \<\Xi(z)-\Xi(x)\>_t^n \le c|z-x|^{2n\alpha}t^{n-1} \int_0^{t} 
  \<G^{\lambda^*}(z-\cdot)+G^{\lambda^*}(x-\cdot)|, \mu_u\>^{2n}du.
 \eeqlb 

Putting together (\ref{mutintegrabilityG}), (\ref{qvntimediff2}), 
Corollary 2.11 of Ethier and Kurtz \cite[Chapter 2]{EthierKurtz86}) and 
Proposition 10.3 of Ethier and Kurtz \cite[Chapter 3]{EthierKurtz86}), 
observe that there is a version of $\{\Xi_t(x):t\ge s\}$ which is a 
time continuous, square-integrable ${\cal F}_t$-martingale for every choice of $x\in {\R}^d$, 
hence we can use the classical martingale inequalities on this version. 

The Burkholder-Davis-Gundy Inequality (Theorem IV.4.1 of Revuz and Yor \cite{RY94}), applied 
to both time continuous, square-integrable ${\cal F}_t$-martingales 
$\{\Xi_t(x)-\Xi_s(x):t\ge s\}$ and $\{\Xi_t(z)-\Xi_t(x):t\ge0\}$, 
yields the existence of a universal constant $C_n > 0$ for each $n\ge1$ 
such that, for every $x,z\in {\R}^d$ and every $0\le s\le t\le T$, there holds 
 \beqlb \lab{BDG}
 \E_{\mu_0} \left( \left| \Xi_t(z)-\Xi_s(x)\right|^{2n} \right)  
 & \leq & 2^{2n-1} \left[ \E_{\mu_0} \left( \left| \Xi_t(z)-\Xi_t(x)\right|^{2n} \right)
  +  \E_{\mu_0} \left( \left| \Xi_t(x)-\Xi_s(x)\right|^{2n} \right) \right]  \nonumber  \\
 & \leq & C_n \left[  \E_{\mu_0} \left( \<\Xi(z)-\Xi(x)\>_t^{n} \right) 
  + \E_{\mu_0} \left( \<\Xi_\cdot(x)-\Xi_s(x)\>_t^{n} \right) \right] \\
 & \leq & c\cdot [|z-x|^{2n\alpha}t^n + (t-s)^n], \nonumber  
 \eeqlb
 this last inequality using (\ref{qvntimediff2}), (\ref{qvnspacediff2}) and Lemma \ref{le_mutintegrability}, 
 incorporating its bounds and $C_n$ into the new $c=c(T, d, h, q^1, n, \alpha, \lambda,\mu_0)>0$ 
 which does not depend on $x$, $z$, $s$ or $t$.
 \qed 

 \bth \lab{th_HolderBranching} 
 Assume that Hypotheses $\ref{hyp:basicassumpFilter}$, $\ref{hyp:basicassumpElliptic}$,  
 $\ref{hyp:basicassumpGaussUniform}$ and $\ref{hyp:basicassumpUniformInteg}$ 
 are satisfied for some $a\ge0$ and $\mu_0\in M_a(\R^d)$.  For either $d=1$, $2$ or $3$, the random field 
 \beqlb \lab{DI1}
 Y_t(x):= \int_0^{t} \int_{\R^d} Q^{\lambda}(x-y)M^\mu(dy,ds)
 \eeqlb
is a square-integrable ${\cal F}_t$-martingale, for every $\lambda > 0$, with quadratic variation given by
\beqnn
  \<Y(x) \>_t = \gamma\sigma^2 \int_0^t \<[Q^{\lambda}(x-\cdot)]^2, \mu_s\> ds 
\eeqnn
 and satisfying $\sup_{x\in\R^d} \E_{\mu_0}  \<Y(x)\>_t  < \infty $ for every $t > 0$.
 There exists a version of $\{Y_t(x), t \geq 0, x \in \R^d \}$ which is 
jointly H\"{o}lder continuous in $(t,x) \in [0, \infty) \times \R^d $, for every $\lambda > 0$.
  \eth

  \proof 
 Since (\ref{DI1}) is from Dawson et al. \cite{DVW19}, we focus once again on the 
 joint H\"{o}lder continuity. Part (ii) of Proposition \ref{prop_mu0integrability} ensures time 
 continuity of stochastic integral $Y_t(x)$ for each fixed $x$. 
 Using H\"{o}lder's inequality directly here, still with conjugate 
 exponents $n$ and $n/(n-1)$ when $n>1$, we get the upper bounds 
 \beqlb \lab{qvntimediff3}
 \<Y_\cdot(x)-Y_s(x)\>_t^{n} \leq(\gamma\sigma^2)^n(t-s)^{n-1}
  \int_{s}^{t} \< [Q^{\lambda}(x-\cdot)]^2, \mu_u\>^{n}du  
 \eeqlb
and  
  \beqlb \lab{qvnspacediff3}
  \<Y(z)-Y(x)\>_t^n \leq(\gamma\sigma^2)^nt^{n-1}\int_0^{t} 
  \< |Q^{\lambda}(z-\cdot)-Q^{\lambda}(x-\cdot)|^2, \mu_u\>^{n}du 
 \eeqlb
 with equality when $n=1$ in both cases. 
Just as in the proof of Theorem \ref{th_HolderEnvironment}, there are constants 
$c=c(\alpha) > 0$ and $c_0=c_0(\alpha) > 0$ for each $\alpha\in(0,1)$, 
verifying (\ref{LSU}) --- with $a_1=c$, $a_2=c_0$ and $\lambda^*=2c_0\lambda$ --- and (\ref{formulaI}), 
such that both 
$Q^{\lambda}(x) \leq c\cdot \tilde{Q}^{\lambda^*}(x)$ and 
$|Q^{\lambda}(z) - Q^{\lambda}(x)|^2 
\leq c|z - x|^{2\alpha}\left[(\tilde{Q}^{\lambda^*}(x))^2+(\tilde{Q}^{\lambda^*}(z))^2\right]$ 
are satisfied, for every $x,z\in {\R}^d$, with $\tilde{Q}^{\lambda^*}$ from (\ref{eqn:greenSBM}).
The rest of the argument follows the proof of Theorem \ref{th_HolderEnvironment}. 
 \qed 
 \\ 
 
 \proof [Proof of Theorem \ref{lt_th1}] 
 All we need to show here is that every term on the right hand side of (\ref{TanakaI}) 
 satisfies individually the H\"{o}lder-type upper bound (\ref{JointHolder}) for their respective norm, as defined 
 by the left hand side. The last two terms do, as a direct consequence of Theorems \ref{th_HolderEnvironment}  
 and \ref{th_HolderBranching}. The other three are less difficult and handled similarly. 
 
Using (\ref{CauchySchwarzHolder}), with $F(x,z,w)= |Q^{\lambda}(z-w)-Q^{\lambda}(x-w)|$ and 
$H(y,w)=I_B(y)$ for any set $B$ with $\lambda_0(B)=1$, the third term (ignoring the constant $\lambda$) 
has its spatial increments bounded, by way of the Cauchy-Schwarz inequality
\[
\left(\int_{s}^{t} \< |Q^{\lambda}(z-\cdot)-Q^{\lambda}(x-\cdot)|, \mu_u\>du\right)^{2n} 
\le (t-s)^{n}\left(\int_{s}^{t} \< |Q^{\lambda}(z-\cdot)-Q^{\lambda}(x-\cdot)|, \mu_u\>^{2}du\right)^{n} , 
\] 
and then (\ref{CauchySchwarzHolder}), to get  
\[
\left(\int_{s}^{t} \< |Q^{\lambda}(z-\cdot)-Q^{\lambda}(x-\cdot)|, \mu_u\>du\right)^{2n} 
\le (t-s)^{2n-1}\int_{s}^{t} \< |Q^{\lambda}(z-\cdot)-Q^{\lambda}(x-\cdot)|, \mu_u\>^{2n}du 
\] 
and the expectation $\E_{\mu_0}$ of the spatial difference taken on both sides is treated identically to (\ref{qvnspacediff3}) 
by setting $s=0$. 

The time difference is treated in the same fashion, using $F(x,z,w)= Q^{\lambda}(x-w)$ and (\ref{qvntimediff3}) instead. 

The first two terms of (\ref{TanakaI}) follow by way of (\ref{LSU}), (\ref{formulaI}) and (\ref{mutintegrabilityQ}), 
just as in the proof of Theorem \ref{th_HolderBranching}. 
\qed 


\section{Proofs}\label{app:ProofsOfLemmas}
\setcounter{equation}{0}

 \subsection{Equivalence of (\ref{nonBAUnif}) and (\ref{nonBA})}\label{app:pf_nonBAUnif}
 \proof 
  We begin by recalling Dawson et al. \cite[Equation 6.2]{DVW19} which states 
 the existence of a constant $C=C(a,d,T) > 0$, for each $a\ge0$, $d\ge1$ and $T > 0$, such that there holds
 \beqlb \lab{mainconstant}
\sup_{0 \leq s \leq T} \sup_{x\in\R^d}\left(I_a^{-1}(x)\int_{\R^d}I_a (y) \varphi_s(y-x) dy\right)\le C.
 \eeqlb 
On any interval $[\epsilon,t]\subset(0,T]$ we have, for all $y\in\R^d$, 
\bea \label{kernelgrowth} 
\sup_{\epsilon \leq s \leq t} \varphi_s(y) \leq \left( \frac{t}{\epsilon} \right)^{d/2} \varphi_t(y) .
\eeqlb
Because of $\lim_{|x|\rightarrow\infty}\varphi_t(x)I_a^{-1}(x)=0$ for any $a$, 
the radius $w=w(d,a,t)>0$ of the smallest ball centered at $0$ 
outside of which $\varphi_t\le I_a$ holds, is such that, for all $x\in\R^d$, we have 
\[
\varphi_t(x) \leq I_a(x)1_{(w,\infty)}(|x|) + (2\pi t)^{-d/2}(1+w^2)^{a/2}I_a(x)1_{[0,w]}(|x|)
\]
where the indicator function of set $N\in{\cal B}(\R)$ is $1_N(x)=1$ if $x\in N$ and $0$ elsewhere. 
Hence there is a constant $c=c(d,a,t)>0$ such that $\varphi_t(x)\le cI_a(x)$ holds everywhere 
and we get  
\[ 
\sup_{y\in\R^d} \sup_{\epsilon \leq s \leq t}\<\varphi_s(y-\cdot) , \mu_0\> 
\leq c \left( \frac{t}{\epsilon} \right)^{d/2} \sup_{y\in\R^d} \<I_a(y-\cdot) , \mu_0\> < \infty
\]  
under Hypothesis \ref{hyp:basicassumpUniformInteg}, using (\ref{kernelgrowth}). 
This is the required bound ensuring that the two statements are equivalent. 
\qed 

 \noindent  {\bf Remark:} 
The need for Hypothesis \ref{hyp:basicassumpUniformInteg} is illustrated next. 
By Fubini's theorem combined with (\ref{kernelgrowth}), any measure $\mu_0\in M_a(\R^d)$ satisfies 
\[
\int_{\R^d}dy\cdot I_a (y) \sup_{\epsilon \leq s \leq t}\<\varphi_s(y-\cdot) , \mu_0\> 
\leq C\left( \frac{t}{\epsilon} \right)^{d/2} \<I_a,\mu_0\> < \infty ,
\] 
hence $ \sup_{\epsilon \leq s \leq t}\<\varphi_s(y-\cdot) , \mu_0\> < \infty$ 
holds $\lambda_0$-almost everywhere in $y\in\R^d$ (and hence on a dense subset of $\R^d$). 
Since the mapping $y\mapsto\varphi_s(y-x)$ is concave for each fixed $x\in\R^d$ and $s>0$, 
it sends any convex combination $y=\alpha w+(1-\alpha)z$ with $\alpha\in[0,1]$ 
to a value satisfying 
$ \min[\varphi_s(w-x), \varphi_s(z-x)]\le\varphi_s(y-x)\le \max[\varphi_s(w-x), \varphi_s(z-x)]$. 
This means that the dense set is also a convex set. Therefore the dense set is all of $\R^d$. 
 
We get, over any ball $B\in{\cal B}(\R^d)$, using the dyadic grid in order to trap 
this ball inside the convex hull of at most $2^d$ dyadic points, 
\[ 
\sup_{y\in B} \sup_{\epsilon \leq s \leq t}\<\varphi_s(y-\cdot) , \mu_0\> < \infty
\]  
for any $\mu_0\in\cup_{a\ge0}M_a(\R^d)$ and under no further assumptions. 
The need for an additional hypothesis in order to control the escape of mass at $\infty$ under $\mu_0$ is now clear. \\ 


 \subsection{Proof of Theorem \ref{tha}}\label{app:pf_th_tha}
 \proof 
The first statement is from Dawson et al. \cite[Theorem 4.3 and Corollary 4.4]{DVW19}.
In the special case of the heat kernel $q_t^{1}(\cdot,y)=\varphi_t(y-\cdot)$ from 
(\ref{gaussiankernel}), a well-known argument using the Fourier transform states that, 
 for any $p\in[1,\infty]$ and $\phi \in L^p(\lambda_0)$, $P^1_t\phi \in C_b^\infty(\R^d)$ holds 
 at all positive times $t > 0$ ---  for instance, see page 180 of Lieb and Loss \cite{LiebLoss01}.
Therefore, under Hypothesis \ref{hyp:basicassumpElliptic}, 
equation (\ref{LSU}) ensures that, for any $m\ge1$, $p\in[1,\infty]$ and 
$Y_0 \in L^p(\lambda_0^m)$, $P^m_t Y_0 \in B((\R^d)^{m})$ holds at all positive times $t > 0$ 
in the general (uniformly elliptic) case as well. Equation (\ref{formulaI}) implies that 
$P^m_{s+t} Y_0=P^m_s P^m_t Y_0 \in C_0((\R^d)^{m})$ also holds 
for all positive times $s > 0$ and $t > 0$, so that even when the starting point of 
process $Y$ is in one of the spaces  $L^p(\lambda_0^m)$, it immediately drops down 
into $C_0((\R^d)^{m})$ and stays there forever $\P_{\mu_0}$-almost surely. 
The comments immediately after (\ref{eqn:Semigroup}) imply that the same is true 
if we replace $C_0((\R^d)^{m})$ by $C_b^2((\R^d)^{m})$,
except at those (finitely many) random jump times $\{\tau_k: 0\le k\le J_0-1\}$, 
which is enough to get bounds in the mean $\E_{\mu_0}$. This is done similarly, by way of 
$P^m_{r+s+t}Y_0=P^m_r P^m_{s+t} Y_0 \in C_b^2((\R^d)^{m})$ which holds 
for all positive times $r > 0$, $s > 0$ and $t > 0$. By Lemma 4.1 of Dawson et al. \cite{DVW19},
there also holds $P^m_{\ell+r+s+t}Y_0=P^m_{\ell} P^m_{r+s+t} Y_0 \in  L^p(\nu_0^m)$, 
for all positive times $\ell > 0$, all $p\in[1,\infty]$ and any measure $\nu_0^m$ satisfying 
Hypothesis \ref{hyp:basicassumpGaussUniform}, thus including both $\lambda_0^m$ and $\mu_0^m$. 
This completes the proof of the second statement, since the $\P_{\mu_0}$-probability of $t$ lying inside 
random set $\{\tau_k: 0\le k\le J_0-1\}$ is $0$.  Inequality (\ref{uniformbound}) now ensues from 
Dawson et al. \cite[Lemma 4.2]{DVW19}, in view of the above, since the Markov property allows us to write 
   \beqnn 
 \| {\cal I}_{a, J_t}^{-1} Y_{t} \|_\infty 
 = \lim_{\epsilon\rightarrow 0} \| {\cal I}_{a, J_t}^{-1} Y_{t-\epsilon}(P^m_\epsilon f)\|_\infty 
 \le c\cdot \lim_{\epsilon\rightarrow 0}\| {\cal I}_{a, m}^{-1} (P^m_\epsilon f)\|_\infty
 \le c\cdot \| {\cal I}_{a, m}^{-1} f\|_\infty 
 \eeqnn
holds $\P_{\mu_0}$-almost surely. Finally, to prove inequality (\ref{Lpmubound}), Jensen's inequality first yields 
$|P^m_t f|^p(x)=|E_x f(X_t)|^p\le E_x |f|^p(X_t)=P^m_t |f|^p(x)$
 for the diffusion $X_t$ associated with semigroup $P^m_t$ and any choice of $f\in L^p(\lambda_0^m)$, 
 $x\in (\R^d)^{m}$, $p\in[1,\infty)$ and $m\ge1$. 
 Since $|\Phi_{ij}^m f|^p(x)=\Phi_{ij}^m |f|^p(x)$ also holds everywhere, we can assume without 
 loss of generality that $p=1$. Inequality (\ref{Lpmubound}) now follows directly from 
 inequality (\ref{uniformbound}).
  \qed
 

 \subsection{Proof of Proposition \ref{prop_mu0integrability}}\label{app:pf_prop_mu0integrability}
 \proof 
   With constants $a_1 > 0$ and $a_2 > 0$ from (\ref{LSU}), we write 
  $S=S_p=(s_1,s_2,\ldots,s_p)$, omitting the index $p$ whenever no ambiguity arises, 
  measure $dS=\prod_{i=1}^{p}ds_i$ and functions $\prod_S=\prod_{i=1}^{p}s_i$, $\Sigma_S=\sum_{i=1}^{p}s_i$
 and $\gamma(S) = a_2\sum_{i=1}^{p} s_i^{-1}$, which satisfies $\gamma(S)\ge a_2p^2/\Sigma_S$.
For any choice of $\mu_0\in M_a(\R^d)$, (\ref{LSU}) yields 
  \beqlb \lab{pnorm}
  (Q^{\lambda}(x))^p \leq \int_{[0,\infty)^p}dS 
  e^{-\lambda\Sigma_S}  \frac{a_1^p}{\prod_S^{d/2}} \exp\{ - |x|^2\gamma(S)\}.    
  \eeqlb

The size of the integrand near $0$ affects the finiteness of this integral significantly so 
we partition $[0,\infty)^p$ into blocks close to and away from the origin. 
For any subset of indices 
$A_n\subset\{1,2,\ldots,p\}$ of cardinal $n\le p$, with $A_0=\emptyset$ 
and complement $A_n^c=\{1,2,\ldots,p\}\setminus A_n$, write $\A_n$ for the 
subset of those points $(s_1,s_2,\ldots,s_p)\in[0,\infty)^p$ such that $s_i\le \eta$
holds for all $i\in A_n$ and $\B_n$ for the one such that 
$s_i >  \eta$ for all $i\in A_n^c$, with $\eta\ge0$ arbitrary but fixed. The symmetry of mappings 
$\prod_S$, $\Sigma_S$ and $\gamma(S)$ in their coordinates makes the integrals to follow, over these sets,  
 dependent only on the cardinality of the sets involved and not on their actual selection, so this 
 ambiguous notation is not problematic and its simplicity fits our purpose. Indeed we can write, 
 for any permutation invariant function of the form 
 $F(S)=\prod_{i=1}^pG(s_i)=\prod_{i\in A_n} G(s_i) \prod_{j\in A_n^c} G(s_j)$ 
 for some $G$ (which is positive in all our applications so let $G$ map $[0,\infty)$ into itself here), 
 \beqlb \lab{inclusionexclusion}
    \int_{[0,\infty)^p} F(S) dS
  & \le &  \int_{\B_0} F(S) dS 
  + \sum_{n=1}^{p-1} \binom{p}{n}  \int_{\A_n\cap\B_n} F(S) dS 
  +  \int_{\A_p} F(S) dS \nonumber \\
  & = & \sum_{n=0}^{p} \binom{p}{n} \left[ \int_{\A_n^*} \prod_{i=1}^{n} G(s_i) ds_i \right] 
            \left[ \int_{\B_n^*} \prod_{j=n+1}^{p} G(s_j) ds_j \right]
  \eeqlb
 where $\binom{p}{n}$ denotes the usual binomial coefficient, 
 $\A_n^*=[0,\eta]^n$ denotes the projection of $\A_n$ onto its constrained coordinates 
 and similarly for $\B_n^*=(\eta,\infty)^{p-n}$. (The integrals over $\A_0^*=\emptyset$ and  
 $\B_p^*=\emptyset$ are set equal to $1$, here and in the rest of the proof.) 
 
 The next step consists in selecting a specific value of $\eta$, as small as possible, in order 
 to control the tail asymptotics of the multiple integral (\ref{pnorm}) decomposed along the lines 
 of (\ref{inclusionexclusion}), that is, its value over the largest sets $\B_n^*$ possible, for every 
 $n\in\{1,2,\ldots,p\}$, thus enabling sharper bounds on as small a neighbourhood of the origin 
 as possible, as this is where the singularity makes integrability most difficult, in the worst case 
 over the set $\A_p^*$. The value $\eta=2a_2p/a$ turns out to be 
 the smallest one for which the following calculations carry through. 
 
So we bound the term in $\B_n^*$ first. 
For every $v\in[0,\infty)$, $\gamma > 0$ and $\beta\ge0$, there holds 
  \beqlb \lab{gamma}
 e^{-\gamma v} (1+v)^{\beta} 
\le 1_{\gamma>\beta} + 1_{\gamma\le\beta} \cdot e^{\gamma} \left( \frac{\beta}{e\gamma} \right)^{\beta} 
\le 1_{\gamma>\beta} + 1_{\gamma\le\beta} \cdot \left( \frac{\beta}{\gamma} \right)^{\beta} .
 \eeqlb
 There ensues successively, since on $\B_n^*$ we have $\min_{\{i=1,2,\ldots,p-n\}} s_i > \eta$ which implies 
 $\gamma(S_{p-n})<a_2(p-n)/\eta=a(p-n)/2p \le \beta=a/2$ when $0 \le n < p$ and $\eta=2a_2p/a$, that there holds, 
 for any $a>0$ and any integers $d\ge1$ and $p\ge1$, 
 \beqlb \lab{tail}
  & & I_a^{-1}(x) \int_{\B_n^*}dS e^{-\lambda\Sigma_S}  \frac{1}{\prod_S^{d/2}}
  \exp\{ - |x|^2\gamma(S)\} \nonumber \\
  & \leq & \int_{\B_n^*}dS e^{-\lambda\Sigma_S}  \frac{1}{\prod_S^{d/2}} 
  \sup_{x\in \R^d} I_a^{-1}(x)\exp\{ - |x|^2\gamma(S)\} \nonumber \\
  & \leq & \int_{\B_n^*}dS e^{-\lambda\Sigma_S}  \frac{1}{\prod_S^{d/2}} 
   \left( \frac{a}{2\gamma(S)} \right)^{a/2} \nonumber \\
  & \leq & \int_{\B_n^*} dS e^{-\lambda\Sigma_S}  \frac{1}{\prod_S^{d/2}} 
  \left( \frac{a\Sigma_S}{2a_2(p-n)^2} \right)^{a/2} \nonumber \\
  & \leq & \frac{1}{p^{d/2}}\frac{1}{(p-n)^a}\left( \frac{a}{2a_2} \right)^{(a+d)/2}  
    \int_{[0,\infty)^{p-n}} dS e^{-\lambda\Sigma_S} \Sigma_S^{a/2} \nonumber \\
  & \leq &\frac{1}{p^{d/2}}\frac{1}{(p-n)^a} \left( \frac{a}{2a_2} \right)^{(a+d)/2}  
    \Gamma(p-n+a/2) \lambda^{-(p-n+a/2)} < \infty ,
  \eeqlb
keeping in mind that for this calculation we write 
  $S=S_{p-n}$, $dS=dS_{p-n}=\prod_{i=1}^{p-n}ds_i$, 
  $\prod_S=\prod_{S_{p-n}}=\prod_{i=1}^{p-n}s_i$, $\Sigma_S=\Sigma_{S_{p-n}}=\sum_{i=1}^{p-n}s_i$
 and $\gamma(S) =\gamma(S_{p-n}) = a_2\sum_{i=1}^{p-n} s_i^{-1}$. 
 Here we write the usual Gamma function $\Gamma (t) = \int_0^\infty v^{t-1}e^{-v}dv$. 
 
Let us denote the last line of (\ref{tail}) by $K(a_2,d,p,n,a,\lambda)$, 
emphasizing that this bound is uniform in $x$ 
in spite of the compensator term $I_a^{-1}(x)$ up front which is required 
due to the generality of starting measure $\mu_0$. 
Define constant $C_p=a_1^p\max_{0 \le n < p}\binom{p}{n}K(a_2,d,p,n,a,\lambda)$, 
eliminate $e^{-\lambda\Sigma_S}\le1$ and 
combine (\ref{pnorm}), (\ref{inclusionexclusion}) and (\ref{tail}) to get, for any time $u\ge0$, 
  \beqlb \lab{pnorm2}
      \<(Q^{\lambda})^p, \mu_u\> 
  & \leq & C_p \sum_{n=0}^{p-1} \int_{\R^d} \mu_u(dx) I_a(x) 
  \left[ \int_{\A_n^*}dS_{n} \frac{1}{\prod_{S_{n}}^{d/2}} \exp\{ - |x|^2\gamma(S_{n})\} \right] \nonumber \\
    & + & a_1^p \int_{\R^d} \mu_u(dx) 
  \left[ \int_{\A_p^*}dS_{p} \frac{1}{\prod_{S_{p}}^{d/2}} \exp\{ - |x|^2\gamma(S_{p})\} \right].  
  \eeqlb
Note the absence of mollifier $I_a(x)\le1$ when $n=p$, making this last term more difficult to control. 
We also have, for every $u\ge0$ and $n\le p$,  
  \beqnn 
\int_{\R^d} \mu_u(dx) \exp\{ - |x|^2\gamma(S_n)\}  
  \le\Upsilon_{\eta^*}(\varphi, \mu_u)\left( \frac{\pi}{\gamma(S_n)} \right)^{d/2} 
  \le\Upsilon_{\eta^*}(\varphi, \mu_u)\left(\frac{\pi}{a_2n}\right)^{d/2}\left(\hbox{$\prod$}{}_{S_n}\right)^{d/2n}
   \eeqnn
where $\Upsilon_{\eta^*}(\varphi, \mu_u)$ is from (\ref{nonBAUnif}) and we used both 
$1/[2\gamma(S_n)]\le\eta/2a_2n\le\eta^*:=\eta/2a_2$ on the set $\A_n^*$ 
and the inequality between the geometric and the harmonic means. 
The double integrals in (\ref{pnorm2}) become, again for every $u\ge0$ and $n\le p$,  
 \beqlb \lab{pnorm3}
  \int_{\R^d} \mu_u(dx) 
  \left[ \int_{\A_n^*}d{S_n} \frac{1}{\prod_{S_n}^{d/2}} \exp\{ - |x|^2\gamma({S_n})\} \right] 
  \leq\Upsilon_{\eta^*}(\varphi, \mu_u)\left(\frac{\pi}{a_2n}\right)^{d/2} \left[ \int_{0}^\eta s^{d/2n-d/2}ds \right]^n 
   \eeqlb
which is finite at $u=0$ under Hypothesis \ref{hyp:basicassumpGaussUniform} 
for all dimensions $d\ge1$ provided $n=1$, for dimensions $d\le3$ when $n=2$
and for all $n\ge1$ when $d=1$ or $2$.

We have just proved that for every measure $\mu_0\in M_a(\R^d)$ satisfying 
Hypothesis \ref{hyp:basicassumpGaussUniform}, there holds 
$\<Q^{\lambda}, \mu_0\> < \infty$ in all dimensions $d\ge1$,  
$\<(Q^{\lambda})^2, \mu_0\> < \infty$ when $d\le3$ and 
$\<(Q^{\lambda})^p, \mu_0\> < \infty$ for all $p\ge1$ when $d=1$ or $2$.  
Notice that Hypothesis \ref{hyp:basicassumpUniformInteg} is not needed for these statements.

Furthermore, recalling that Theorem \ref{tha} ensures that $\mu_u\in M_a(\R^d)$ 
holds $\P_{\mu_0}$-almost surely for any $u\ge0$ as soon as $\mu_0\in M_a(\R^d)$, 
the remarks following Hypothesis (\ref{hyp:basicassumpGaussUniform}) imply that 
$\mu_u(w-\cdot)\in M_a(\R^d)$ then also holds $\P_{\mu_0}$-almost surely,  
for every $w\in\R^d$ simultaneously. 
Substituting it for $\mu_u$ in the calculations starting with (\ref{pnorm2}), 
ensures that the upper bound in (\ref{pnorm3}) is valid uniformly in $w$, 
since the shift invariance in (\ref{nonBAUnif}) entails 
\[
\Upsilon_{\eta^*}(\varphi, \mu_u(w-\cdot))=\Upsilon_{\eta^*}(\varphi(w-\cdot), \mu_u)
=\Upsilon_{\eta^*}(\varphi, \mu_u).
\] 
We can therefore also write 
$\sup_{w\in \R^d} \<Q^{\lambda}(w-\cdot), \mu_0\> < \infty$ in all dimensions $d\ge1$,  
$\sup_{w\in \R^d} \<[Q^{\lambda}(w-\cdot)]^2, \mu_0\> < \infty$ when $d\le3$ and 
$\sup_{w\in \R^d} \<[Q^{\lambda}(w-\cdot)]^p, \mu_0\> < \infty$ for all $p\ge1$ when $d=1$ or $2$, 
provided Hypothesis \ref{hyp:basicassumpUniformInteg} attains.
 
To get $\sup_{w\in \R^d} \int_{\R^d}|\partial_{x_i}Q^{\lambda}(w-x)|  \mu_0(dx) < \infty$ for all $d\ge1$,
 set values $p=1$, $n=0$ and $\eta=2a_2p/a$ in (\ref{tail}) with the added term $s^{-1/2}$ to obtain 
\[
  I_a^{-1}(x) \int_{\eta}^\infty ds e^{-\lambda s} s^{-(d+1)/2} \exp\{ - a_2|x|^2/s\}
  \leq K(a_2,d,1,0,a+1,\lambda) < \infty ,
\]
valid for all $x\in\R^d$, with the constant $K(a_2,d,1,0,a+1,\lambda)>0$ defined right after (\ref{tail}). 
Next, successively using (\ref{eqn:interchange}), (\ref{LSU}), $\eta^*:=\eta/2a_2$, 
$C_1=a_1K(a_2,d,1,0,a+1,\lambda)>0$, Hypothesis \ref{hyp:basicassumpGaussUniform} 
and the same line of reasoning as above implies finiteness at $u=0$ of    
  \beqlb \lab{partialboundat0}
  \int_{\R^d}|\partial_{x_i}Q^{\lambda}(w-x)|  \mu_u(dx)  
  & \leq & \int_0^{\infty} e^{-\lambda s} \frac{a_1}{s^{(d+1)/2}} \<\exp\{ - a_2|\cdot|^2/s\}, \mu_u(w-\cdot)\> ds 
     \nonumber \\
  & \leq & C_1 \<I_a(w-\cdot), \mu_u\>
+ a_1\Upsilon_{\eta^*}(\varphi, \mu_u)\left( \frac{\pi}{a_2} \right)^{d/2} \Gamma(\frac{1}{2}).
  \eeqlb
 
When combined with Theorem \ref{tha} all finiteness statements above propagate in time, that is, 
they remain true in the mean $\E_{\mu_0}$ if we replace $\mu_0$ by $\mu_t$. 

Indeed, duality identity (\ref{dual}) simplifies to $\E_{\mu_0} \<f, \mu_t\> = \<P^{1}_{t}f, \mu_0\>$ 
for any $f\in L^1(\mu_0)$. Choosing such an $f\ge0$ which also satisfies $f\in L^1(\lambda_0)$,
there holds 
\[
\<P^{1}_{t}f(w-\cdot), \mu_0\>\le\Upsilon_t(q^{1}(w-\cdot), \mu_0)\<f, \lambda_0\>
=\Upsilon_t(q^{1}, \mu_0(w-\cdot))\<f, \lambda_0\>\le A^*\Upsilon_{ct}(\varphi, \mu_0)\<f, \lambda_0\>
\] 
 using (\ref{Aronsonbounds}) under Hypothesis \ref{hyp:basicassumpGaussUniform} and  
 the other finiteness statements follow through the appropriate selection for $f$, under 
 Hypothesis \ref{hyp:basicassumpUniformInteg}.
  \qed


\bibliography{refers}
\bibliographystyle{amsplain}

\end{document}